\documentclass{sn-jnl}

\makeatletter
\let\cl@chapter\undefined
\makeatletter

\usepackage{float}
\usepackage{amsmath,amsfonts,amssymb}
\usepackage[usenames,dvipsnames]{xcolor}
\usepackage{url}
\usepackage{booktabs}
\usepackage{hyperref}
\usepackage{caption, subcaption}
\usepackage{cleveref}
\usepackage[short]{optidef}
\usepackage{suffix}
\usepackage{algorithm, algorithmicx}
\usepackage[noEnd]{algpseudocodex}
\usepackage{pgfplots}
\usepackage{comment}
\usepackage{longtable}
\usepackage{lscape}

\newtheorem{theorem}{Theorem}[section]
\theoremstyle{definition}
\newtheorem{definition}[theorem]{Definition}
\newtheorem{example}[theorem]{Example}
\theoremstyle{remark}
\newtheorem{remark}{Remark}

\newcommand{\R}{{\mathbb{R}}}

\newcommand{\I}{{\mathcal{I}}}

\newcommand{\A}{{\mathcal{A}}}

\newcommand{\Lag}{{\mathcal{L}}}

\newcommand{\BNLP}{{\mathrm{BNLP}}}

\newcommand{\zcc}{\zeta}
\newcommand{\zcco}{\zeta_1}
\newcommand{\zcct}{\zeta_2}

\newcommand{\epstol}{\epsilon_{\mathrm{tol}}}
\newcommand{\pnorm}[2]{\norm{#2}_{#1}}
\newcommand{\infnorm}[1]{\pnorm{\infty}{#1}}
\newcommand{\twonorm}[1]{\pnorm{2}{#1}}

\newcommand{\ncc}{n_{\mathrm{cc}}}
\newcommand{\nzero}{n_0}

\newcommand{\apr}{\alpha_{\mathrm{pr}}}
\newcommand{\adu}{\alpha_{\mathrm{du}}}
\newcommand{\aprmax}{\alpha_{\mathrm{pr}}^{\textrm{max}}}
\newcommand{\adumax}{\alpha_{\mathrm{du}}^{\textrm{max}}}

\newcommand{\transp}[1]{#1^{\top}}
\newcommand{\inv}[1]{#1^{-1}}

\newcommand{\LPEC}{{\mathrm{LPEC}}}
\newcommand{\PROJLPEC}{{\mathrm{PROJ-LPEC}}}
\newcommand{\LPCC}{{LPCC}}
\newcommand{\LPCCs}{{LPCCs}}
\newcommand{\MPEC}{{MPEC}}
\newcommand{\MPECs}{{MPECs}}
\newcommand{\MPCC}{{MPCC}}
\newcommand{\MPCCs}{{MPCCs}}

\newcommand{\ipopt}{\texttt{IPOPT}}
\newcommand{\madnlp}{\texttt{MadNLP}}
\newcommand{\knitro}{\texttt{Knitro}}
\newcommand{\solver}{\texttt{CCOpt}}

\newcommand{\casadi}{\texttt{CasADi}}
\newcommand{\gurobi}{\texttt{Gurobi}}
\newcommand{\highs}{\texttt{HiGHS}}

\newcommand{\lcqpow}{\texttt{LCQPow}}
\newcommand{\MacMPEC}{\texttt{MacMPEC}}
\newcommand{\NOSBENCH}{\texttt{NOSBENCH}}
\newcommand{\nosnoc}{\texttt{nosnoc}}

\DeclareMathOperator{\diag}{diag}

\newcommand{\paren}[1]{\left( #1 \right)}
\newcommand{\brac}[1]{\left[ #1 \right]}
\newcommand{\cbrac}[1]{\left\{ #1 \right\}}
\WithSuffix\newcommand\paren*[1]{( #1 )}
\WithSuffix\newcommand\brac*[1]{[ #1 ]}
\WithSuffix\newcommand\cbrac*[1]{\{ #1 \}}
\newcommand{\norm}[1]{{\left\Vert#1\right\Vert}}
\newcommand{\Set}[2]{\left\{\, #1 \mid #2\,\right\}}
\WithSuffix\newcommand\Set*[2]{\{\, #1 \mid #2\,\}}

\makeatletter
\newcommand\footnoteref[1]{\protected@xdef\@thefnmark{\ref{#1}}\@footnotemark}
\makeatother

\newcommand{\Input}{\Statex \textbf{Input:}}
\newcommand{\Desc}[2]{\Statex \hspace{\algorithmicindent} \makebox[7em][l]{#1}\parbox[t]{0.8\textwidth}{#2}}

\title{{\solver}: an Open-Source Solver for Large-Scale Mathematical Programs with Complementarity Constraints}
\author*[1]{\fnm{Anton} \sur{Pozharskiy}}\email{anton.pozharskiy@imtek.uni-freiburg.de}
\author[2]{\fnm{Fran\c{c}ois} \sur{Pacaud}}\email{francois.pacaud@minesparis.psl.eu}
\author[1,3]{\fnm{Moritz} \sur{Diehl}}\email{moritz.diehl@imtek.uni-freiburg.de}
\author[1]{\fnm{Armin} \sur{Nurkanovi\'c}}\email{armin.nurkanovic@imtek.uni-freiburg.de}

\affil[1]{\orgdiv{Department of Microsystems Engineering (IMTEK)}, \orgname{University of Freiburg}, \orgaddress{\city{Freiburg}, \country{Germany}}}
\affil[2]{\orgdiv{Centre Automatique et Syst\`emes}, \orgname{Mines Paris-PSL}, \orgaddress{\city{Paris}, \country{France}}}
\affil[3]{\orgdiv{Department of Mathematics}, \orgname{University of Freiburg}, \orgaddress{\city{Freiburg}, \country{Germany}}}

\date{\today}

\begin{document}
\abstract{
 This paper presents the Julia package {\solver}, built on top of the interior-point solver {\madnlp}.
 {\solver} implements a suite of algorithms for Mathematical Programs with Complementarity Constraints (\MPCCs).
 The solver additionally comes with interfaces for use in Matlab, Python, and C++.
 {\MPCCs} have recently gained renewed attention in engineering optimization, as complementarity provides a powerful modeling tool for nonsmooth functions and logical conditions.
 These problems are inherently challenging since their nonlinear programming reformulations violate classical regularity conditions at all feasible points, complicating both theoretical analysis and numerical treatment.
 Consequently, specialized algorithms are required to handle this degeneracy, and several approaches have been proposed.
 We implement a toolbox of methods, including relaxation and penalty approaches, as well as a crossover to recently proposed active-set methods.
 Our solver is based on nonlinear interior-point algorithms that couple the relaxation or penalty parameter with the barrier parameter, yielding substantial speedups compared to standard implementations.
 Both monotone and nonmonotone strategies for updating this joint parameter update are proposed and investigated.
 In addition, we propose regularization techniques that improve the conditioning of the KKT system for small relaxation parameters, enhancing robustness and computational efficiency.
The implementation is validated on the classical \texttt{MacMPEC} benchmark, large-scale problems in security-constrained optimal power flow, optimal control of nonsmooth systems, as well as on quadratic programs with complementarity constraints arising in model predictive control.
This benchmarking reveals an algorithmically driven improvement of often an entire order of magnitude over other methods, including commercial solvers.
}
\keywords{large-scale numerical optimization, mathematical programs with complementarity constraints (MPCCs), interior-point method, open-source software, nonlinear programming}
\maketitle

\section{Introduction}
\label{sec:intro}
This paper regards Mathematical Programs with Complementarity Constraints\\ (\MPCCs) of the following form:
\begin{mini!}[2]
  {\substack{x\in\R^n}}
  {f(x)}
  {\label{eq:mpec}}
  {}
  \addConstraint{c(x)}{=0\label{eq:mpec_eq}}
  \addConstraint{0}{\le x_0\label{eq:mpec_lb0}}
  \addConstraint{0\le x_1}{\perp x_2 \ge 0,\label{eq:mpec_comp}}
\end{mini!}
with the partition of variables $x=(x_0,x_1,x_2) \in \R^n$, with $x_0 \in \R^{\nzero}$, $x_1,x_2\in \R^{\ncc}$.
The functions $f:\R^n \to \R$ and $c:\R^n \to \R^{m} $ are assumed to be at least twice continuously differentiable.
The notation $0 \leq x_1 \perp x_2 \geq 0$ means that all components for the vectors $x_1$ and $x_2$ must be non-negative, and that vectors $x_1$ and $x_2$ are orthogonal, i.e., either $x_{1,i} = 0$ or $x_{2,i} = 0$ for all $i$.
For ease of exposition we do not include general inequality constraints in~\eqref{eq:mpec}.
If inequality constraints are present, and defined via at least twice continuously differentiable functions, all results are readily extended via the addition of slack variables to $x_0$.

The complementarity constraints \eqref{eq:mpec_comp} make the theoretical and computational aspects of the optimization problem~\eqref{eq:mpec} more difficult than in standard nonlinear programming~\cite{Scheel2000}.
If the constraints \eqref{eq:mpec_comp} are replaced by a set of inequality constraints, then a standard Nonlinear Program (NLP) is obtained:
\begin{mini!}
  {\substack{x\in\R^n}}
  {f(x)}
  {\label{eq:mpec_nlp}}
  {}
  \addConstraint{c(x)}{=0\label{eq:mpec_nlp_ineq}}
  \addConstraint{0}{\le x_0\label{eq:mpec_nlp_lb0}}
  \addConstraint{x_{1,i}}{\geq 0,\ x_{2,i} \geq 0,\ x_{1,i} x_{2,i} \leq 0, \label{eq:mpec_nlp_comp}}{\ i = 1,\ldots,\ncc}
\end{mini!}
However, this NLP violates standard constraint qualifications such as the \sloppy Mangasarian-Fromovitz constraint qualification (MFCQ) at all feasible points~\cite{Jane2005,Scheel2000}.
Consequently, the set of its Lagrange multipliers is unbounded~\cite{Gauvin1977}.
As a practical consequence, standard NLP methods applied to~\eqref{eq:mpec_nlp} become inefficient~\cite{Kim2020,Nurkanovic2024b}.
The lack of regularity also renders the application of the classical Karush–Kuhn–Tucker (KKT) conditions to~\eqref{eq:mpec_nlp} problematic, complicating both the definition and verification of stationarity~\cite{Scheel2000}.
Several alternative reformulations of the complementarity constraints~\eqref{eq:mpec_comp} have been proposed, including those based on C-functions~\cite{Facchinei2003,Leyffer2006a} or equality constraints for the bilinear term in \Cref{eq:mpec_nlp_comp}.
However, none of these approaches eliminates the degeneracy inherent in the NLP reformulation~\eqref{eq:mpec_nlp}.

{\MPCCs} occur naturally in a vast variety of scientific and engineering disciplines,
including robotics~\cite{Wensing2023}, process engineering~\cite{Baumrucker2008,Raghunathan2003}, optimal control of nonsmooth systems~\cite{Nurkanovic2023f,Nurkanovic2022a}, power systems~\cite{Pacaud2025}, sparsity optimization~\cite{Feng2018,Kanzow2024}, and bi-level optimization problems~\cite{Kim2020}.
Further, complementarity constraints provide a generic framework for reformulating many nonsmooth constraints and objectives~\cite{Hegerhorst2020,Szmuk2020,Raghunathan2003}.
Many of these problems are often solved ``ad-hoc'' by practitioners.
With {\solver} we seek to provide a robust and performant solver which can be used to accelerate these applications significantly.

\paragraph{Related work:}
In order to solve \MPCCs\ and tackle the difficulty posed by the degeneracy of complementarity constraints, several approaches have been proposed.
These can be categorized in two major classes of algorithm: \emph{regularization} and \emph{active set} methods.

Regularization methods operate on the notion of replacing the offending complementarity constraints with a better behaved formulation depending on a regularization parameter.
A sequence of regularized NLPs is solved to recover the solution of the original \MPCC\ at the limit.
These methods can be split into three subclasses: relaxation, smoothing, and penalization methods.
\emph{Relaxation methods} replace the constraint $x_{1,i}x_{2,i} \le 0$ (and sometimes the corresponding lower bounds) with the inequality $\Phi_{\mathcal{I}}(x_{1,i}, x_{2,i}; \tau) \le 0$, which corresponds to the L-shaped complementarity set when $\tau=0$.
The goal is to transform the \MPCC\ into a sequence of regularized NLPs which satisfy constraint qualifications (often MFCQ) for all $\tau > 0$ and solve a sequence of such problems with $\tau\rightarrow 0$.
The oldest of these methods, proposed by S. Scholtes in 2001, is the so called Scholtes' global relaxation: $\Phi_{\mathcal{I}}(a,b;\tau) = ab -\tau$~\cite{Scholtes2001}.
Since then, many more relaxation schemes have been proposed with a variety of convergence properties~\cite{Kadrani2009,Kanzow2013,Lin2003,Steffensen2010,DeMiguel2005,Wang2023}.
\emph{Smoothing methods} are similar and replace the constraint in \Cref{eq:mpec_nlp_comp}, with a C-function $\Phi_{\mathcal{E}}(x_{1,i}, x_{2,i}; \tau) = 0$ encoding the L-shaped complementarity set when $\tau=0$.
These have been shown to have generally worse accuracy characteristics when compared to inequality based relaxation methods, however~\cite{Ralph2004}.
\emph{Penalty methods}, on the other hand, move the bilinear constraints into the objective with a penalty parameter $\rho$.
The simplest version penalizes the complementarity residual $\rho \transp{x_1}x_2$ and maintains the lower bounds on these variables~\cite{Anitescu2005,Ralph2004,Leyffer2006}.
Extensions to these methods have been proposed with more advanced penalty formulations with the aim to tame the nonconvexity of the penalty term~\cite{Lin2024b,Wang2025}.
A numerical comparison and survey of these methods can be found in~\cite{Nurkanovic2024b,Hoheisel2013,Kanzow2015,Kim2020}.

Active set methods try to take advantage of the combinatorial nature of \MPCCs\ and find solutions by selecting the active branches of the set of complementarity constraints.
They then approximately solve a sequence of well-posed NLPs based on a sequence of complementarity active sets to find a stationary point.
Examples of such algorithms can be found in~\cite{Nurkanovic2025,Kazi2024,Izmailov2008,Jiang1999}.
In addition, direct application of the Sequential Quadratic Programming (SQP) method to~\eqref{eq:mpec_nlp} has been studied in~\cite{Fletcher2006}.

The present work is closely related to the algorithm proposed by Raghunathan and Biegler in~\cite{Raghunathan2005} which uses an interior-point method in tandem with the Scholtes' global relaxation approach.
This algorithm drives the relaxation and barrier parameter to zero together, solving the \MPCC\ via a sequence of approximate NLP solutions.
Unfortunately, no implementation of these methods remains available presently.
Our new algorithm improves on~\cite{Raghunathan2005} by allowing both the barrier and relaxation parameters to be adjusted independently and in a nonmonotonous fashion.
Further, we implement several specialized regularizations of the KKT system which tackle the indefinite terms coming from the nonconvex relaxed complementarity constraints.
Along with this we improve upon the Vicente-Wright~\cite{Vicente2002} regularization scheme used in~\cite{Raghunathan2005} to address the degeneracy of the Jacobian, via a novel bound relaxation scheme.
Finally, we provide a performant implementation of this algorithm based on the state-of-the-art interior-point framework {\madnlp}, and provide interfaces to the modeling framework CasADi~\cite{Andersson2019} via a C interface, creating a fully open source software environment.

Another closely related set of works are~\cite{DeMiguel2005,Wang2023} which present algorithms using two-sided relaxations of the complementarity feasible set.
These algorithms relax both the lower bounds and the upper bounds a la Scholtes' global relaxation and chose which of the relaxations to drive to zero depending on estimated multipliers.
We take a similar approach but avoid certain degenerate cases which can often occur by allowing iterates to escape the non-negative quadrant $(x_1,x_2) \geq 0$.

\begin{remark}[A comparison to non-interior methods]
    Our approach differs from the ``non-interior'' point method used classically to solve problems with complementarity constraints~\cite{Chen1993,Lin2022}.
    In those methods, both the upper level and lower level complementarities of the KKT conditions are reformulated using a so-called C-function, such as the well known Fischer-Burmeister function or its many modifications~\cite{Chen2000a}.
    C-functions have zero level sets which correspond exactly to the L-shaped complementarity set.
    Unlike interior-point methods,
    these C-functions do not require keeping $(x,z)$ in the nonnegative quadrant and therefore avoid the need for the ``fraction-to-boundary'' rule which limits step size in interior-point methods.
    However, they reformulate the complementarity constraints with equality constraints,
    decreasing the accuracy achieved at a given relaxation parameter $\tau$ to $O(\tau^{\frac{1}{4}})$,
    while inequality relaxations yield $O(\tau^{\frac{1}{2}})$~\cite{Ralph2004}.
    This theoretical observation has been corroborated with our previous numerical experiments~\cite{Pozharskiy2023}.
    Further, we are unaware of any performant implementation of non-interior-point methods for {\MPCCs}.
\end{remark}

\paragraph{Contributions:}
The following details the main contributions of this paper:
\begin{enumerate}
\item In \Cref{sec:madnlpc} we describe a novel interior-point algorithm for solving \MPCCs.
  This algorithm is based on the basic form proposed by Raghunathan and Biegler in~\cite{Raghunathan2005}, and drives the barrier parameter $\mu$ and global Scholtes complementarity relaxation parameter $\tau$ to zero jointly.
\item In \Cref{sec:madnlpc:hessian_regularization} we introduce two novel methods for handling the structural indefiniteness of the Lagrangian Hessian which occurs due to the bilinear constraint terms.
  These methods take advantage of the complementarity structure and are used to minimize the number of matrix factorizations needed in each iteration of the algorithm.
\item In \Cref{sec:madnlpc:endgame} we introduce a novel lower bound relaxation scheme which aims to alleviate the degeneracy in the constraint Jacobian which occurs as the complementarity relaxation parameter approaches zero.
  This relaxation scheme uses estimates of the {\MPCC} multipliers to adjust bounds on the complementarity variables $x_1$ and $x_2$ conservatively in order to prevent iterates from escaping the non-negative orthant.
\item In \Cref{sec:madnlpc:homotopy_update} we introduce several novel parameter update rules for our algorithm and discuss their impacts on the algorithm.
  The first of these aims to solve the more relaxed NLPs to higher accuracy which empirically leads to finding lower objective solutions often with a reduced number of iterations required to converge.
  The second and third describe MPCC-tailored extensions to the adaptive barrier update schemes implemented in the software package \texttt{LOQO}~\cite{Shanno2000} and {\ipopt}~\cite{Waechter2006}.
    \item In parallel to the Scholtes relaxation, \solver\ also supports penalty regularizations.
    In particular, in \Cref{sec:penalty}, we describe our implementation of an interior penalty method similar to the algorithm described by Leyffer et al. in~\cite{Leyffer2006}.
    We include a version of our novel Lagrangian Hessian regularization in \Cref{sec:ell1:hess_regularization}.
    \item \Cref{sec:benchmarks} proceeds to validate the sizable, often order-of-magnitude, performance improvements that the implementation achieves over other methods for solving \MPCCs.
    The comparison includes problems from several major fields of interest, including large-scale problems from optimal control, huge-scale problems arising from security-constrained optimal power flow, as well as quadratic programs with complementarity constraints coming from real-time algorithms for model predictive control.
    {\solver} compares favorably not only to other open-source solvers but even to state-of-the-art commercial solvers, including those with specialized algorithms for {\MPCC}.
    \item To allow easy access to our high performance solver for a large section of the scientific and engineering community we provide not only a Julia package but also a compiled shared library interface
    which can be called from any language that provides a foreign function interface.
    In addition, we implement an interface for \solver\ to the popular modeling and automatic differentiation framework CasADi~\cite{Andersson2019}, enabling the use from Matlab, Python, and C++.
\end{enumerate}

\paragraph{Outline:}
We structure the rest of this paper as follows.
\Cref{sec:background} is dedicated to basic background on {\MPCCs}.
In this section we also review the underlying interior-point architecture of our implementation.
\Cref{sec:madnlpc}, \Cref{sec:penalty}, and \Cref{sec:crossover}, describe the algorithms implemented in {\solver}, including an interior-point method, an interior-penalty method, as well as a crossover to an {\MPCC} active-set method, respectively.
We also compare some of the options available in our implementation and provide some guidance to practitioners.
We then proceed in \Cref{sec:benchmarks} to benchmark our implementation against classical approaches on several suites of challenging test problems.
Finally, in \Cref{sec:conclusion}, we discuss conclusions and lay out future research directions.

\paragraph{Notation:}
Subscripts are used to identify components of vectors or matrices, e.g., $x_i$, and superscripts to denote iterates, e.g., $x^k$.
Functions that are evaluated at particular iterates are denoted as $f^k := f(x^k)$.
The concatenation of two vectors $a \in \R^n$ and $b \in \R^m$ is shortly written as $(a,b) = \begin{bmatrix} a^\top & b^\top\end{bmatrix}^\top$.
The same notation is adapted accordingly to the concatenation of several vectors.
Given two scalars $a,b\in\R$, $\max(a,b)$ returns the larger of them.
Let $a, b \in \R^n$, then $\max(a)$ returns the largest component of $a$, and $c=\max(a,b)$ is a vector in $\R^n$ with the components $c_i = \max(a_i,b_i),\ i =1,\ldots,n$.
The $\min$ function is defined analogously.
Let $x\in\R^{n}$, then $\diag(x)\in\R^{n\times n}$ returns a diagonal matrix with $x_i$ as its diagonal elements.
We define the set $\R_{\ge 0}^n\coloneq\Set*{x\in\R^n}{x\ge 0}$.
A vector $e$ is used to represent the vector of ones of the appropriate size, $I$ is the identity matrix of appropriate size, and where necessary $0$ is the zero matrix of appropriate size.
The Jacobian of a function is denoted by $\nabla c(x)\coloneq \frac{\partial c}{\partial x}(x)^\top$.


\section{Background}
\label{sec:background}
In this section, we recall the first-order stationarity conditions tailored to {\MPCCs}.
Then we first briefly review the interior-point method for nonlinear programming.

\subsection{Optimality conditions for {\MPCCs}}
First, we define some index sets for the complementarity constraint of the \MPCC~\eqref{eq:mpec}.
Denote the feasible set of the {\MPCC}~\eqref{eq:mpec} by
\[
\Omega = \{ x\in \R^n \mid x_0 \ge 0, \; c(x) = 0, \; 0 \leq x_1 \perp x_2 \geq 0\}.
\]
We regard the following index sets, defined at a feasible point $x$:
\begin{align*}
    \I_{+0}(x) &=   \{i \in \{1,\ldots,\ncc\} \mid x_{1,i}>0, x_{2,i}=0\},\\
    \I_{0+}(x) &= \{i \in \{1,\ldots,\ncc\}\mid x_{1,i}=0, x_{2,i}>0\},\\
    \I_{00}(x) &= \{i \in \{1,\ldots,\ncc\}\mid x_{1,i}=0, x_{2,i}=0\}.
\end{align*}
Most of the computational and theoretical difficulties arise if the \textit{biactive set} $\I_{00}(x)$ is nonempty.

Regard the tight NLP (TNLP), a standard NLP, derived for the {\MPCC}~\eqref{eq:mpec} at a feasible point $x^*\in\Omega$:
\begin{mini!}
  {\substack{x \in \R^{n}}}
  {f(x)\label{eq:tnlp_objective}}
  {\label{eq:tnlp}}
  {}
  \addConstraint{c(x)}{=0\label{eq:tnlp_ineq}}
  \addConstraint{x_0}{\ge 0}
  \addConstraint{x_{1,i}}{= 0,\ x_{2,i} \geq 0,\qquad\qquad\label{eq:tnlp_branch_h}}{\forall i \in  \mathcal{I}_{0+}(x^*)}
  \addConstraint{x_{1,i}}{\ge 0,\ x_{2,i} = 0,\qquad\qquad\label{eq:tnlp_branch_g}}{\forall i \in   \mathcal{I}_{+0}(x^*)}
  \addConstraint{x_{1,i}}{= 0,\ x_{2,i} = 0,\qquad\qquad\label{eq:tnlp_biactive}}{\forall i \in     \mathcal{I}_{00}(x^*).}
\end{mini!}
If $x^*\in \Omega$ is a local minimizer of the {\MPCC}~\eqref{eq:mpec}, then it is a local minimizer of the TNLP~\cite{Scheel2000}.
If \eqref{eq:tnlp_biactive} is replaced by $x_{1,i} \geq 0,\ x_{2,i} \geq 0,\; \forall  i \in   \mathcal{I}_{00}(x^*)$, then we obtain the so-called relaxed NLP (RNLP).
Let us recall some {\MPCC}-tailored notions.

\begin{definition}\label{def:mpec_cqs}
    The {\MPCC} \eqref{eq:mpec} is said to satisfy the {\MPCC} Linear Independence Constraint Qualification, {\MPCC}-LICQ, (respectively {\MPCC}-MFCQ) at a feasible point $x^*\in \Omega$ if the corresponding $\mathrm{TNLP}$ \eqref{eq:tnlp} satisfies the LICQ (respectively MFCQ) at the same point $x^*$.
\end{definition}
Next, we define the \MPCC\ Lagrangian, which corresponds to the standard Lagrangian of the TNLP:
\begin{align}\label{eq:mpcc-lag}
    \Lag^{\mathrm{\MPCC}}(x,y,z,\zcc) = f(x) +\transp{y}c(x)- \transp{z}_0 x_0- \transp{\zcc}_1 x_1 - \transp{\zcc}_2 x_2,
\end{align}
where $y\in\R^{m}$ $z_0 \in \R_{\ge 0}^{\nzero}$, $\zcco \in \R^{\ncc}$ and $\zcct \in \R^{\ncc}$ are the {\MPCC} Lagrange multipliers.

By applying the KKT conditions to the TNLP~\eqref{eq:tnlp}, one obtains the W-stationarity concept, and applying them to the RNLP yields S-stationary~\cite{Scheel2000}.
Further restrictions on the multipliers
$ \zcc^*_{1,i}, \zcc^*_{2,i}, \ i \in \I_{00}(x^*)$, lead to M-, C-, and A-stationarity, which are not directly related to the KKT conditions of an auxiliary NLPs, see e.g.~\cite{Kim2020,Outrata1999,Scheel2000}.
These concepts are summarized in the next definition.
\begin{definition}[Stationarity concepts for {\MPCCs}]\label{def:mpec_stationarity}
    For a feasible point $x^* \in \Omega$, we distinguish the following stationarity concepts.
    \begin{itemize}
        \item Weak stationarity (W-stationarity) \cite{Scheel2000}:
        A point $x^*$ is called W-stationary if the corresponding TNLP~\eqref{eq:tnlp} admits the satisfaction of the KKT conditions, i.e., there exist Lagrange multipliers $\mu^*,\nu^*$ and $\xi^*$ such that:
        \begin{align*}
            &\nabla_x \mathcal{L}^{\mathrm{\MPCC}}(x^*,y^*,z^*, \zcc^*) = 0,\\
            &0 \leq z_0^* \perp x_0 \geq 0,\\
            & x^{*}_{1,i} \geq0, \zcc^*_{1,i} = 0, \; &\forall i \in    \I_{+0}(x^*),\\
            & x^{*}_{2,i} \geq0, \zcc^*_{2,i} = 0,\; &\forall i \in     \I_{0+}(x^*),\\
            & x^{*}_{1,i} =0,\; \zcc_{1,i}^* \in \R, \; &\forall i \in  \I_{0+}(x^*)\cup\I_{00}(x^*),\\
            &x^{*}_{2,i} =0,\; \zcc_{2,i}^* \in \R, \; &\forall i \in   \I_{+0}(x^*)\cup\I_{00}(x^*).
        \end{align*}
        \item Strong stationarity (S-stationarity) \cite{Scheel2000}:
        A point $x^*$ is called S-stationary if it is weakly stationary and $\zcc_{1,i}^* \geq 0, \zcc_{2,i}^* \geq 0$ for all $i \in \I_{00}(x^*)$.
        \item Clarke stationarity (C-stationarity) \cite{Scheel2000}:
        A point $x^*$ is called C-stationary if it is weakly stationary and $\zcc_{1,i}^*\zcc_{2,i}^* \geq 0$ for all $i \in \I_{00}(x^*)$.
        \item Mordukhovich stationarity (M-stationarity) \cite{Outrata1999}:
        A point $x^*$ is called \\M-stationary if it is weakly stationary and if either $\zcc_{1,i}^* >0$ and $\zcc_{2,i}^* >0$ or $\zcc_{1,i}^*\zcc_{2,i}^* =0$ for all $i \in \I_{00}(x^*)$.
        \item Abadie stationarity (A-stationarity) \cite{Flegel2005a}:
        A point $x^*$ is called A-stationary if it is weakly stationary and $\zcc_{1,i}^* \geq 0$ or $\zcc_{2,i}^* \geq0$ for all $i \in \I_{00}(x^*)$.
    \end{itemize}
\end{definition}
Note that in the case of $\I_{00}(x^*) = \emptyset$, then all stationarity concepts reduce to S-stationarity, which is the most restrictive, since it has the tightest condition in the multipliers $\zcc_{1,i}^*, \zcc_{2,i}^*, i \in \I_{00}(x^*)$.
Weaker stationarity concepts, in descending order of strength, are M, C, A, and W, with the stronger ones implying the weaker ones.

A generic first-order optimality condition, independent of the regularity of the NLP, can be stated with use of the tangent cone, which is defined as follows.
The tangent cone at $x\in\Omega$ to the set $\Omega$ is defined as $\mathcal{T}_{\Omega}(x)
= \{ d\in \R^n \mid \exists \{x^k\} \subset \Omega, \{t^k\}\subset \R_{\geq0}: \lim\limits_{k\to \infty} t^k = 0, \lim\limits_{k\to \infty} x^k = x, \lim\limits_{k\to \infty} \frac{x^k-x}{t^k} = d\}$.
More specifically a first-order necessary optimality conditions via the tangent cone for \eqref{eq:mpec} are given by the following theorem~\cite[Corollary 3.3.1]{Luo1996}.
\begin{theorem}\label{th:b_stationarity}
  Let $x^*\in\Omega$ be a local minimizer of \eqref{eq:mpec}, then it holds that
  \begin{align}\label{eq:geometric_b_stationarity}
    \nabla f(x^*)^\top d \geq 0\; \textrm{ for all } d \in \mathcal{T}_{\Omega}(x^*),
  \end{align}
  or equivalently, $d = 0$ is a local minimizer of the following optimization problem:
  \begin{align}\label{eq:b_stationariry}
    \min_{d \in \R^{n}} \quad &   \nabla f(x^*)^\top d\quad     \textnormal{s.t.} \quad   d \in \mathcal{T}_{\Omega}(x^*).
  \end{align}
\end{theorem}
If a point $x^*$ satisfies the condition above, it is said that ``geometric'' B-stationarity holds~\cite{Flegel2005a,Luo1996,Scheel2000}.
If the so-called {\MPCC}-Abadie constraint qualification is satisfied~\cite{Flegel2005a}, then we can write an algebraic condition for B-stationarity:
\begin{definition}[B-stationarity~\cite{Flegel2005a,Luo1996}]\label{def:b_stationarity}
  A point $x^*$ is called algebraically B-stationary if $d=0$ solves the following {\LPCC}:
  \begin{mini!}[2]
    {\substack{d \in \R^n}}
    {\nabla f(x^*)^\top d}
    {\protect\label{eq:lpec_reduced_theory}}
    {}
    \addConstraint{c(x^*)+ \nabla c(x^*)^\top d}{\ge 0,\quad}{\forall i \in \A(x^*)}
    \addConstraint{x_{1,i}^* + d_{1,i} }{=0,\quad}{\forall i \in  \I_{0+}(x^*)}
    \addConstraint{x_{2,i}^* + d_{2,i}}{=0,\quad}{\forall i \in  \I_{+0}(x^*)}
    \addConstraint{0 \leq  x_{1,i}^* + d_{1,i}  \perp x_{2,i}^* + d_{2,i}}{\ge 0,\quad}{\forall i \in \I_{00}(x^*).}
  \end{mini!}
\end{definition}

The multiplier-based stationarity concepts and B-stationarity are related as follows.
If $x^*$ is an S-stationary point of the \sloppy {\MPCC}~\eqref{eq:mpec}, then it is also B-stationary.  If in addition, the {\MPCC}-LICQ holds at $x^*$, then every B-stationary point is S-stationary~\cite[Theorem 4]{Scheel2000}.
Notably, if {\MPCC}-LICQ does not hold, there might exist B-stationary points that are not S-stationary.
This assumption cannot be relaxed, as the next weaker condition, {\MPCC}-MFCQ is already not sufficient for S-stationarity, for a counterexample see \cite[Example 3]{Scheel2000}.

\subsection{An interior-point skeleton}\label{sec:ipm}
The two major algorithms in {\solver} both constitute extensions of a filter line-search Interior Point Method (IPM) for constrained nonlinear programming (NLP)~\cite{Waechter2006,Shin2024}.
In this section we discuss the skeleton of such an algorithm, the pseudo-code for which can be found in \Cref{alg:ipm}.

The general problem solved by IPMs for NLP takes the form:
\begin{mini!}
  {\substack{x\in\R^{n}}}
  {f(x)}
  {\protect\label{eq:ip}}
  {}
  \addConstraint{c(x)}{= 0}
  \addConstraint{x}{\ge 0.}
\end{mini!}
We assume that the NLP (\ref{eq:ip}) is regular, in the sense it satisfies the Mangasarian-Fromowitz Constraint Qualification (MFCQ) and the Second Order Sufficient Condition (SOSC) at a solution $x^*$.
These conditions are commonly required by IPM algorithms as they are necessary and sufficient for the local convergence theory developed in the literature to hold~\cite{Vicente2002}.
Any NLP with general inequality constraints can be reformulated into \eqref{eq:ip} without loss of generality by adding nonnegative slacks to any inequality constraints.
The Lagrangian of (\ref{eq:ip}) is:
\begin{equation*}
  \Lag(x,y,z) = f(x) + \transp{y}c(x) - \transp{z}x,
\end{equation*}
where $y\in\R^m$ and $z\in\R^n$ are the Lagrange multipliers of the equality constraints and bound multipliers respectively.
The Karush-Kuhn-Tucker (KKT) conditions for the NLP~\eqref{eq:ip} read as:
\begin{align*}
  \nabla_x \Lag(x,y,z) &= 0,\\
  c(x) &= 0,\\
  0\le x \perp z &\ge 0.
\end{align*}
Interior point methods then approximately solve a sequence of smooth perturbed KKT systems of the form:
\begin{equation}
  \label{eq:nlp_stationarity}
  F_{\mu}(w) = \begin{bmatrix}
    \nabla_x\Lag(x,y,z)\\
    c(x)\\
    Xz - \mu^k
  \end{bmatrix} = 0,
\end{equation}
with $\mu^k$ being a homotopy barrier parameter which is driven to zero, and ${X = \diag(x)}$, while maintaining $(x,z) >0$.
A generic globalized two phase IPM algorithm for (\ref{eq:ip}) can be decomposed into five main steps:
\begin{enumerate}
\item Termination criteria (Line 2 of \Cref{alg:ipm}) for evaluating whether the current point $x^k$ is a stationary point of the given NLP, up to some tolerance.
\item An update rule for the relaxation homotopy (Line 3 of \Cref{alg:ipm}).
For NLPs this corresponds to an update rule for $\mu^k$, the relaxation parameter for the complementarity slackness found in the perturbed KKT conditions in \Cref{eq:nlp_stationarity}.
\item Step calculation (Line 4 of \Cref{alg:ipm}), which is classically done via evaluating a Newton step in the (perturbed) first-order stationarity conditions of the problem in \Cref{eq:nlp_stationarity}, though other methods also exist.
\item A globalization scheme (Line 5 of \Cref{alg:ipm}) to ensure the method globally converges to the minimizer of the NLP.
\item A recovery scheme (Line 7 of \Cref{alg:ipm}) which attempts to recover from points $x^k$ where step calculation and globalization fail to deliver an acceptable step.
\end{enumerate}

For the final two steps we use the filter line-search method and restoration phase respectively, as implemented in \ipopt~\cite{Waechter2009} and \madnlp~\cite{Shin2024}.
We set a maximum fraction-to-boundary step size via the formulas:
\begin{subequations}
  \label{eq:frac_to_boundary}
  \begin{align}
    \aprmax &= \max\Set{\alpha\in[0,1]}{x^k + \alpha d^k \ge (1-\eta_b)x^k}\\
    \adumax &= \max\Set{\alpha\in[0,1]}{z^k + \alpha z^k \ge (1-\eta_b)z^k}
  \end{align}
\end{subequations}
for an $\eta_b \in (0,1)$ to maintain $(x^k, z^k) > 0$.
This prevents taking steps too close to the boundaries too quickly.
In the following sections we describe modifications to steps one through three to specialize our method to handling complementarity constraints.
\begin{algorithm}[t]
  \caption{An Interior Point Skeleton.}
  \label{alg:ipm}
  \begin{algorithmic}[1]
    \State $k = 0$, $x^0$, $y^0$, $z^0$ initialization
    \While{Termination criteria are not satisfied.}
      \State Update homotopy parameters.
      \State Evaluate step direction $d^k=(d_x^k,d_y^k,d_z^k)$.
      \State Select primal and dual step size $\apr$ and $\adu$ by some globalization scheme.
      \If{Step evaluation or step size selection fails}
        \State Find $x^{k+1}$, permitting progress, via recovery scheme cf. the ``restoration phase''~\cite{Waechter2006}.
      \Else
        \State $x^{k+1} = x^k + \apr d_x^k$, $y^{k+1} = y^k + \adu d_y^k$, $z^{k+1} = z^k + \adu d_z^k$
      \EndIf
      \State $k = k+1$
    \EndWhile
  \end{algorithmic}
\end{algorithm}

\section{A relaxation-based interior-point algorithm for MPCCs}
\label{sec:madnlpc}
In this section we describe an interior-point relaxation method, which solves a series of relaxed problems with the relaxation parameter $\tau$ driven to zero in tandem with the barrier parameter $\mu$.
Several algorithmic specifics are then explored in detail, including: specialized Hessian regularization, joint relaxation and barrier parameter updates, as well as an endgame algorithm which relaxes lower bounds in order to alleviate terminal numerical issues.
Finally, the crossover to the active-set algorithm described in~\cite{Nurkanovic2025} is detailed.
This method allows us to ensure convergence to true local minima satisfying B-stationarity.

\subsection{Relaxation-based algorithm}
Inspired by~\cite{Raghunathan2005} we implement a filter line-search interior-point method which solves the relaxed \MPCC\ by updating the relaxation parameter $\tau$ alongside the interior-point barrier parameter $\mu$.
For simplicity of exposition we treat the following NLP:
\begin{mini!}[2]
  {\substack{x\in\R^{n},\ s\in\R^{\ncc}}}
  {f(x)}
  {\protect\label{eq:madnlpc_nlp}}
  {}
  \addConstraint{c(x)}{=0}
  \addConstraint{x}{\ge 0}
  \addConstraint{X_1x_2 + s - \tau e}{=0}
  \addConstraint{s}{\ge 0,}
\end{mini!}
which is the Scholtes relaxation of the NLP (\ref{eq:mpec_nlp}).
Recall that $x$ are the primal decision variables, partitioned into $x = (x_0,x_1,x_2)$, $s$ are the slacks for the Scholtes inequality constraints, and $\tau > 0$ is the relaxation parameter for the inequalities in \Cref{eq:mpec_nlp_comp}.
The function $c:\R^n\rightarrow\R^m$ represents the original nonlinear equality constraints.
We remind the reader that the matrix $X_1 = \mathrm{diag}(x_1)$.
The Lagrangian is then defined as:
\begin{equation}
  \label{eq:madnlpc_lag}
  \Lag(x,s,y,z) = f(x) + \transp{c(x)}y_c + \transp{(X_1 x_2 + s - \tau e)}y_s - \transp{z_0}x_0 - \transp{z_1}x_1 - \transp{z_2}x_2 - \transp{z_s}s,
\end{equation}
in which we introduce multipliers $y_c\in\R^{n_c}$, $y_s\in\R^{\ncc}$, $z_1\in\R^{\ncc}$, $z_2\in\R^{\ncc}$, and $z_s\in\R^{\ncc}$.
Note that $z_1$ and $z_2$ are not the same multipliers as $\zcc_1$ and $\zcc_2$ in \Cref{eq:mpcc-lag}, but the bound multipliers for the corresponding decision variables $x_1$ and $x_2$.
For more compact notation we use $y = (y_c,y_s)$, $z_x = (z_0,z_1,z_2)$, $z=(z_x,z_s)$ and $X_i = \diag(x_i)$, $S = \diag(s)$, and $Z_i = \diag(z_i)$, $W_{ij} = \nabla_{x_i,x_j}(f(x)+\transp{c(x)}y_c)$, and the constraint Jacobian $J_i = \transp{\nabla_{x_i}c(x)}$.
The details are provided in \Cref{alg:madnlpc}, and the rest of this section is dedicated to describing the implementation of its components.
The comments in \Cref{alg:madnlpc} provide pointers to the particular section which discuss the corresponding details.
\begin{algorithm}[t]
  \caption{An interior-point relaxation algorithm for {\MPCCs}.}
  \label{alg:madnlpc}
  \begin{algorithmic}[1]
    \State $k = 0$
    \While{\eqref{eq:term_crit:term} is not satisfied}
      \State $\mu^k = \texttt{update\_mu}(x^k,s^k,y^k,z^k)$
      \Comment{\Cref{sec:madnlpc:homotopy_update}}
      \State $\tau^k = \texttt{update\_tau}(x^k,s^k,y^k,z^k,\mu^k)$
      \Comment{\Cref{sec:madnlpc:homotopy_update}}
      \If{$r(x^k,s^k,y^k,z^k)\le \epsilon_{\mathrm{endgame}}$}
        \State Execute the endgame algorithm.
        \Comment{\Cref{sec:madnlpc:endgame}}
      \EndIf
      \State Evaluate KKT matrix $K^k$ and KKT residual $r^k$.
      \Comment{\Cref{sec:madnlpc:kkt_regularization}}
      \State Calculate search direction $d^k$ by solving $K_a^kd_a^k = r_a^k$ in \cref{eq:aug_kkt}.
      \Comment{\Cref{sec:madnlpc:step_calculation}}
      \State Calculate step lengths $\apr$ and $\adu$ using a filter line-search method.
      \If{step sizes are too small or search direction calculation fails}
        \State Begin restoration phase for fixed $\tau$ a la \cite{Waechter2006}.
      \Else
        \State $x^{k+1} = x^k + \apr d_x^k$, $s^{k+1} = s^k + \apr d_s^k$
        \State $y^{k+1} = y^k + \adu d_y^k$, $z^{k+1} = z^k + \adu d_z^k$
        \State $k=k+1$
      \EndIf
    \EndWhile
  \end{algorithmic}
\end{algorithm}

\subsection{Step calculation}\label{sec:madnlpc:step_calculation}
The KKT conditions for the relaxed MPCC~\eqref{eq:madnlpc_nlp} read as:
\begin{align*}
  \nabla_{x_0}f(x) - \transp{J_0}y_c &= 0,\\
  \nabla_{x_1}f(x) - \transp{J_1}y_c + X_2y_s - z_1&= 0,\\
  \nabla_{x_2}f(x) - \transp{J_2}y_c + X_1y_s - z_2 &= 0,\\
  y_s - z_s &= 0,\\
  c(x) &= 0,\\
  X_1 x_2 + s - \tau e &= 0,\\
  0\le x_0 \perp z_0 &\ge 0,\\
  0\le x_1 \perp z_1 &\ge 0,\\
  0\le x_2 \perp z_2 &\ge 0,\\
  0\le s \perp z_s &\ge 0.
\end{align*}
In an IPM, the complementarity constraints are smoothed using the barrier parameter~$\mu$, resulting in the perturbed KKT conditions:
\begin{subequations}
  \label{eq:ipm_kkt_system}
  \begin{alignat}{2}
    \nabla_{x_0}f(x) + \transp{J_0}y_c &\eqqcolon r_1 = 0,\\
    \nabla_{x_1}f(x) + \transp{J_1}y_c + X_2y_s - z_1 &\eqqcolon r_2 = 0,\\
    \nabla_{x_2}f(x) + \transp{J_2}y_c + X_1y_s - z_2 &\eqqcolon r_3 = 0,\\
    y_s - z_s &\eqqcolon r_4 = 0,\\
    c(x) &\eqqcolon r_5 = 0,\\
    X_1 x_2 + s - \tau e &\eqqcolon r_6 = 0,\label{eq:ipm_kkt_system:scholtes}\\
    X_0z_0 - \mu e &\eqqcolon r_7 = 0,\\
    X_1z_1 - \mu e &\eqqcolon r_8 = 0,\\
    X_2z_2 - \mu e &\eqqcolon r_9 = 0,\\
    Sz_z - \mu e &\eqqcolon r_{10} = 0,\\
    x, s, z &> 0.
  \end{alignat}
\end{subequations}
The interior-point scheme proceeds by taking Newton steps on these equations, which requires solving the linear KKT system:

\begin{equation}
  \label{eq:unreduced_kkt}
  \begin{bmatrix}
    W_{00} & W_{01} & W_{02} & 0 & \transp{J}_0 & 0 & -I & 0 & 0 & 0\\
    W_{00} & W_{11} & W_{12}+Y_s & 0 & \transp{J}_1 & X_2 & 0 & -I & 0 & 0\\
    W_{20} & W_{21}+Y_s & W_{22} & 0 &\transp{J}_2 & X_1 & 0 & 0 & -I & 0\\
    0 & 0 & 0 & 0 & 0 & I & 0 & 0 & 0 & -I\\
    J_0 & J_1 & J_2 & 0 & 0 & 0 & 0 & 0 & 0 & 0\\
    0 & X_2 & X_1 & I & 0 & 0 & 0 & 0 & 0 & 0\\
    Z_0 & 0 & 0 & 0 & 0 & 0 & X_0 & 0 & 0 & 0\\
    0 & Z_1 & 0 & 0 & 0 & 0 & 0 & X_1 & 0 & 0\\
    0 & 0 & Z_2 & 0 & 0 & 0 & 0 & 0 & X_2 & 0\\
    0 & 0 & 0 & Z_s & 0 & 0 & 0 & 0 & 0 & X_s
  \end{bmatrix}
  \begin{bmatrix}
    \Delta x_0\\
    \Delta x_1\\
    \Delta x_2\\
    \Delta s\\
    \Delta y_c\\
    \Delta y_s\\
    \Delta z_0\\
    \Delta z_1\\
    \Delta z_2\\
    \Delta z_s
  \end{bmatrix}
  = -
  \begin{bmatrix}
     r_1\\
     r_2\\
     r_3\\
     r_4\\
     r_5\\
     r_6\\
     r_7\\
     r_8\\
     r_9\\
     r_{10}
  \end{bmatrix}.
\end{equation}
As in most IPM implementations for NLPs, we do not solve this system directly but rather transform it into the so called \textit{Augmented KKT system}~\cite{Nocedal2006}.
This is done by eliminating the rows in \cref{eq:unreduced_kkt} corresponding to the bound multipliers $z$ with the relations:
\begin{alignat*}{2}
  \Delta z_0 &= -\inv{X}_0(r_7 + Z_0\Delta x_0), \quad \Delta z_1 &= -\inv{X}_1(r_8 + Z_1\Delta x_1),\\
  \Delta z_2 &= -\inv{X}_2(r_9 + Z_2\Delta x_2), \quad \Delta z_s &= -\inv{S}(r_{10} + Z_s\Delta s),
\end{alignat*}
that yields the augmented KKT system:
\begin{equation}
  \label{eq:aug_kkt}
  \begin{bmatrix}
    W_{00}+\Sigma_0 & W_{01} & W_{02} & 0 & \transp{J}_0 & 0\\
    W_{00} & W_{11}+\Sigma_1 & W_{12}+Y_s & 0 & \transp{J}_1 & X_2\\
    W_{20} & W_{21}+Y_s & W_{22}+\Sigma_2 & 0 &\transp{J}_2 & X_1\\
    0 & 0 & 0 & \Sigma_s & 0 & I\\
    J_0 & J_1 & J_2 & 0 & 0 & 0\\
    0 & X_2 & X_1 & I & 0 & 0
  \end{bmatrix}
  \begin{bmatrix}
    \Delta x_0\\
    \Delta x_1\\
    \Delta x_2\\
    \Delta s\\
    \Delta y_c\\
    \Delta y_s
  \end{bmatrix}
  = -
  \begin{bmatrix}
     r_1+\inv{X}_0r_7\\
     r_2+\inv{X}_1r_8\\
     r_3+\inv{X}_2r_9\\
     r_4+\inv{S}r_{10}\\
     r_5\\
     r_6
  \end{bmatrix},
\end{equation}
where $\Sigma_i = \inv{X}_iZ_i$ and $\Sigma_s = \inv{S}Z_s$.
We write this in a compact form, separating out the contributions from bound multipliers and the Scholtes relaxation multipliers to the Hessian block:
\begin{equation}
  \label{eq:augmented_kkt_small}
  Kd = -r \; ,
\end{equation}
where:
\begin{equation*}
  K =
  \begin{bmatrix}
    W + Q & 0 &\transp{B}\\
    0 & \Sigma_s & \transp{C}\\
    B & C & 0
  \end{bmatrix},\quad
  d =
  \begin{bmatrix}
    \Delta x\\
    \Delta s\\
    \Delta y
  \end{bmatrix}, \quad
  r=
  \begin{bmatrix}
     r_{x}\\
     r_4+\inv{S}r_{10}\\
     r_y
  \end{bmatrix},
\end{equation*}
\begin{equation*}
  W = \nabla^2_x(f(x)+\transp{c(x)}y_c)=
  \begin{bmatrix}
    W_{00} & W_{01} & W_{02}\\
    W_{00} & W_{11} & W_{12}\\
    W_{20} & W_{21} & W_{22}
  \end{bmatrix},\quad
  Q =
  \begin{bmatrix}
    \Sigma_0 & 0 & 0\\
    0 & \Sigma_1 & Y_s\\
    0 & Y_s & \Sigma_2
  \end{bmatrix},
\end{equation*}
\begin{equation*}
  B =
  \begin{bmatrix}
    J_0& J_1& J_2\\
    0 & X_2 & X_1
  \end{bmatrix},\quad
  C =
  \begin{bmatrix}
    0 \\
    I
  \end{bmatrix},\quad
  r_x =
  \begin{bmatrix}
    r_1+\inv{X}_0r_7\\
    r_2+\inv{X}_1r_8\\
    r_3+\inv{X}_2r_9
  \end{bmatrix},\quad
  r_y =
  \begin{bmatrix}
    r_5\\
    r_6
  \end{bmatrix},
\end{equation*}

\subsection{Termination criteria}\label{sec:madnlpc:term_crit}
To determine the termination of this algorithm we check the following criteria in at each iteration:
\begin{subequations}
  \label{eq:term_crit}
  \begin{align}
    l(x,s,y,z) &=
        \begin{bmatrix}
          \nabla_x\Lag(x,s,y,z)\\
          c(x)\\
          Xz_x\\
          Sz_s\\
          X_1x_2
        \end{bmatrix},\label{eq:term_crit:stat}\\
    \infnorm{l(x,s,y,z)} &\le \epstol,\label{eq:term_crit:term}
  \end{align}
\end{subequations}
where $l$ is the stationarity conditions of the NLP (\ref{eq:madnlpc_nlp}) along with the upper level complementarity measure $X_1x_2$.
As classically done in homotopy based approaches for \MPCCs~\cite{Lin2003,Kadrani2009,Nurkanovic2024b}, we make the choice to evaluate the multiplicative complementarity residual instead of $\max(x_1,x_2)$.
Due to the limitations of double precision floating point calculation, solving these problems to better than a multiplicative residual of $10^{-8}$ (that is $\infnorm{X_1x_2} \le \tau = 10^{-8}$) is impractical.
This, however, would be required if we were to use $\max(x_1,x_2)$ as the complementarity measure in the termination criteria as in the case of biactive constraints $i\in \I_{00}$ it holds that ${\max(x_{1,i},x_{2,i})} = O(\sqrt{\tau})$.
In an attempt to overcome this accuracy barrier, we optionally resort to an active-set crossover strategy discussed in \Cref{sec:crossover}.

\subsection{KKT matrix regularization}
\label{sec:madnlpc:kkt_regularization}
There are several issues which can arise when attempting to solve the augmented KKT system in \cref{eq:aug_kkt}, exacerbated by the structure of the relaxed complementarity constraints.
Firstly, note that the off-diagonal contribution of the multipliers $y_s$ (seen as the $Y_s$ component of $Q$) makes an indefinite contribution to the Hessian block $W+Q$.
While it is known that the off-diagonal elements of $Q$ do not contribute to the reduced Hessian at the solution~\cite{Fletcher2006,DeMiguel2005}, during intermediate iterations the reduced Hessian may become significantly indefinite.
Traditionally, this is handled by an \emph{inertia-correction} step (briefly described in \Cref{sec:madnlpc:inertia_correction}) which may require many expensive refactorizations to identify a sufficient perturbation to recover a positive definite reduced Hessian.
We propose an additional step before falling back to the inertia-correction routine which attempts to rectify the indefinite contribution of the $Q$ matrix, requiring only one refactorization.

The second issue is the ill-conditioning of the KKT system as the solver approaches a solution.
It is well known the diagonal contribution of the active bounds in the augmented Hessian block $Q$ grows to infinity close to the solution~\cite{Shin2024}.
This increases the condition number of the matrix, leading to significant difficulty in computing accurate search directions.
This problem is made worse by the fact that for all of the Scholtes-relaxed constraints the $\Sigma_s$ block grows to $\infty$ as $\tau$ approaches zero.
One insightful way to view this problem is to note that as $\tau$ approaches zero the relative interior of the feasible set begins to disappear and the constraints $x_1x_2 -\tau\le 0$ and (without loss of generality) $x_1\ge 0$ get closer to violating MFCQ.
We propose a lower bound relaxation approach inspired by the findings in~\cite{DeMiguel2005} which attempts to enlarge the interior of the feasible set without compromising the validity of the solution.
The details of this approach are discussed in \Cref{sec:madnlpc:endgame}.

\subsubsection{Inertia correction}\label{sec:madnlpc:inertia_correction}
Note that the linear system described in \Cref{eq:aug_kkt}, may not necessarily be invertible if the reduced Hessian, that is the Hessian projected onto the active constraints, is indefinite.
Further, if the constraint Jacobian is not full rank, which happens structurally in {\MPCCs}, also leads to a failure to calculate a decent direction for the problem.
To combat this we use the implementation of the regularization found in \madnlp~\cite{Shin2024}, and based on prior art in IPM, namely, IPOPT~\cite{Waechter2005,Nocedal2006}, adding scaled identities to the KKT matrix:
\begin{equation}
  \label{eq:aug_kkt_ic}
  \begin{bmatrix}
    W + Q +\delta_w I & 0 &\transp{B}\\
    0 & \Sigma_s + \delta_w I & \transp{C}\\
    B & C & -\delta_cI
  \end{bmatrix}
  \begin{bmatrix}
    \Delta x\\
    \Delta s\\
    \Delta y
  \end{bmatrix}
  =-
  \begin{bmatrix}
     r_{x}\\
     r_4+\inv{S}r_{10}\\
     r_y
  \end{bmatrix},
\end{equation}
The values $\delta_w\ge 0$ and $\delta_c\ge 0$ are selected either by an inertia correction method, or by an inertia free method implemented in \madnlp.
The inertia correction method uses the ``inertia-revealing'' nature of the LBL factorization to evaluate whether the KKT matrix has the correct number of positive and negative eigenvalues (in our case $n+\ncc$ positive eigenvalues and $m+\ncc$ negative eigenvalues).
Verifying the updates of $\delta_w$ and $\delta_c$ via checking the KKT matrix' inertia requires expensive refactorizations and can often be a dominant component of the solver's runtime, particularly for MPCCs.

\subsubsection{Complementarity constraints Hessian regularization}\label{sec:madnlpc:hessian_regularization}
Due to the off diagonal contribution of the Scholtes constraint multipliers, the Hessian of the Lagrangian $W+Q$ may become indefinite.
This may require a large primal regularization $\delta_w$ computed during the inertia-correction step of the IPM.
Inertia-correction is an iterative algorithm requiring a numerical factorization of the KKT matrix for each adjustment of $\delta_w$ and $\delta_c$.
As such, we implement an additional regularization scheme operating specifically on the bottom right 2-by-2 block of the matrix $Q$.
First note that this block can be rearranged into the block diagonal matrix
\begin{equation}
  Q =
  \begin{bmatrix}
    B_1 & 0 & \cdots & 0\\
    0 & B_2 & & \vdots\\
    \vdots & &\ddots &\vdots \\
    0 &\cdots &\cdots & B_{\ncc}
  \end{bmatrix},\quad
  B_i =
  \begin{bmatrix}
    \inv{x}_{1,i}z_{1,i} & y_{s,i}\\
    y_{s,i} & \inv{x}_{2,i}z_{2,i}
  \end{bmatrix}.
\end{equation}
The $2\times 2$ matrices $B_i$ can be regularized in several ways.

\paragraph{Regularization 1 (critical multiplier).}
We note that each block is positive definite if and only if the inequality:
\[
  z_{1,i}z_{2,i}\ge y_{s,i}^2x_{1,i}x_{2,i},
\]
holds.
We can therefore calculate a maximal $\bar{y}_{s,i} \le \sqrt{\frac{z_{1,i}z_{2,i}}{x_{1,i}x_{2,i}}}$, which maintains the positive-definiteness of block $B_i$ (and therefore of $Q$).
When using this method we replace $B_i$ with the matrix
\begin{equation}
  \label{eq:B_bar}
  \bar{B}_i =
  \begin{bmatrix}
     \inv{x}_{1,i}z_{1,i} & \alpha_B\bar{y}_{s,i}\\
    \alpha_B\bar{y}_{s,i} & \inv{x}_{2,i}z_{2,i}
  \end{bmatrix},
\end{equation}
with $\alpha_B = 0.999$, for each $i=1,\ldots,\ncc$.

\paragraph{Regularization 2 (eigenvalue decomposition).}
The second regularization technique operates on each 2-by-2 block $B_i$ individually using an eigenvalue decomposition (which is trivial for 2-by-2 matrices):
\begin{equation}
  \label{eq:B_eig}
  B_i = V_i\Lambda_i\inv{V}_i
\end{equation}
where $\Lambda_i$ is the 2-by-2 diagonal matrix with the eigenvalues of $B_i$ on its diagonal.
We then reassemble $\hat{B}_i$ with a modified $\hat{\Lambda}_i$ whose diagonal elements are clipped to be at minimum of $\lambda_{\mathrm{min}} = 10^{-8}$:
\begin{equation}
  \label{eq:B_hat}
  \hat{B}_i = V_i\hat{\Lambda}_i\inv{V}_i
\end{equation}
We default to using the regularization in \Cref{eq:B_bar} as in practice this approach demonstrates better performance on many problems.
It is observed that the approach using the eigenvalue decomposition will frequently cause large elements to appear on the diagonal of $\hat{B}_i$ which can cause severe numerical issues for step computation, which outweighs the benefits of the regularization scheme.

It is important to note that, in the general case, $Q$ may not be the only source of indefiniteness in the reduced Hessian.
We still rely on the iterative inertia-control algorithm~\cite{Waechter2006} to account for the non-complementarity sources of indefiniteness in $W$, as well as a possibly rank-deficient constraint Jacobian.
As such, we implement the above described Hessian regularization scheme as the first step of the inertia-control algorithm, that is if the factorization of the initial, unmodified, KKT matrix yields an incorrect inertia, we apply the regularization scheme and attempt to refactorize.
If this fails we fall back to the standard inertia-correction scheme.

\begin{remark}
  In the case of well-posed LPCCs and quadratic programs with complementarity constraints (namely with positive semi-definite Hessians, and satisfying MPCC-LICQ) the above regularization of $Q$ is sufficient to ensure that the corresponding KKT matrix has the correct inertia.
  For these problems it is possible to disable inertia-correction entirely: we obtain at most two KKT matrix factorizations per iteration.
\end{remark}

\subsubsection{A lower bound relaxation endgame}\label{sec:madnlpc:endgame}
In the Scholtes inequality relaxation approach, the interior of the feasible set of the problem disappears as $\tau\rightarrow 0$.
This leads to numerical instability in step computation and often leads to a rapid increase in the magnitude of the Lagrange multipliers for the complementarity lower bounds and Scholtes relaxation, as well as numerical issues due to the near LICQ violation.
We observe that these effects can lead to poor terminal convergence behavior in both classical Scholtes homotopy methods and the algorithm proposed in~\cite{Raghunathan2005}.
One approach to mitigating this issue is the Vicente-Wright regularization approach.
We propose an alternative scheme taking advantage of the specific structure of an {\MPCC}.

An interesting observation made by DeMiguel et al in~\cite{DeMiguel2005}, is that if a problem satisfies MPCC Strong Second Order Sufficient Conditions (MPCC-SSOSC) at a solution to the Scholtes relaxed MPCC, at most two of $x_1\ge 0$, $x_1x_2 -\tau \le 0$, and $x_2\ge 0$ can be active, which can be determined from the signs of the complementarity multipliers ($\zcc_1$ and $\zcc_2$ in \Cref{eq:mpcc-lag}).
This insight means that only at most two of the inequality constraints representing the relaxation of $0\le x_1 \perp x_2 \ge 0$ are actually needed at the solution.
To take advantage of this we introduce two more algorithmic relaxation parameters per complementarity constraint: $\delta_1\ge 0$ and $\delta_2\ge 0$, as well as a new Scholtes relaxation parameter for each complementarity: $\tau_v \in\R^{\ncc}$ .
Doing this we write the complementarity reformulation:
\begin{subequations}\label{eq:two_sided}
  \begin{align}
    0&\le x_1 + \delta_1,\\
    0&\le x_2 + \delta_2,\\
    X_1x_2-\tau_v&\le 0.
  \end{align}
\end{subequations}
Note that with $\delta_1 = \delta_2 = 0$ and all elements of $\tau_v$ equal, we recover exactly the Scholtes relaxation formulation.
The algorithm described by DeMiguel et al. then takes some nonzero initial values for $\delta_1$, $\delta_2$, and $\tau_v$ and takes Newton steps with the corresponding KKT system.
At each iteration, for each complementarity, one of  $\delta_1$, $\delta_2$, or $\tau_v$ is driven towards zero based (only) on the current estimate of the complementarity multipliers.

In our experience, this approach fails to deliver correct results when true minima of the relaxed problem exist at points in the interior of the relaxed feasible set and requires some heuristics to escape such points.
Further, the relaxation of variable bounds at a solution can lead to reported solutions arbitrarily far from the real solution to a problem, particularly in the case where the objective is nonconvex or the constraints are nearly linearly dependent.
However, we observe that the effect of a non-disappearing strict interior can improve convergence in the dual variables significantly.
As such, we attempt to mitigate these effects by introducing an ``endgame'' approach where we only relax the bounds via $\delta_1$, and $\delta_2$ when the KKT error (\cref{eq:term_crit}) is sufficiently small and the estimates of the multipliers $\zcco$ and $\zcct$ are sufficiently accurate.
The correspondence of solutions to this lower bound relaxed problem to solutions to the original \MPCC\ follows directly from \cite[Theorem 3.1]{DeMiguel2005}.
The endgame algorithm is described in \Cref{alg:endgame}.
We trigger this endgame algorithm, if requested by the user, when the total KKT residual $r(x,s,y,z)$ is smaller than some threshold (by default $\epsilon_{\mathrm{endgame}} = 10^{-6}$).

\begin{algorithm}[t]
  \caption{An endgame bound relaxation algorithm.}
  \label{alg:endgame}
  \begin{algorithmic}[1]
    \State $\hat{\zcc}_1 = z_1 - z_sx_2$
    \State $\hat{\zcc}_2 = z_2 - z_sx_1$
    \For{$i=1$ to $\ncc$}
    \If{$\hat{\zcc}_{1,i} \le -\Xi(x,y,z)$}
    \Comment{Relax lower bound of $x_{1,i}$}
    \State $\delta_{1,i} = \min(\Psi(\hat{\zcc}_1, x_2,\vec{\tau}_i), \delta_{\textrm{max}}$
    \ElsIf{$\hat{\zcc}_{2,i} \le -\Xi(x,y,z)$}
    \Comment{Relax lower bound of $x_{2,i}$}
    \State $\delta_{2,i} = \min(\Psi(\hat{\zcc}_{2,i}, x_{1,i},\vec{\tau}_i), \delta_{\textrm{max}}$
    \Else
    \State Do nothing to $\delta_1$ and $\delta_2$.
    \EndIf
    \EndFor
  \end{algorithmic}
\end{algorithm}

A key property of the endgame algorithm is it must mitigate the escape of solver iterates from the non-negative orthant.
Additionally, the algorithm should not relax the lower bound unless the estimated {\MPCC} multipliers are of sufficient magnitude.
To this end we introduce two functions $\Xi(x,y,z)$ and $\Psi(\hat{z}, x, \tau)$.
The first function $\Xi(x,y,z)$ is used to evaluate how negative the estimated MPCC multipliers must be before relaxing the corresponding lower bound.
We implement the function:
\begin{equation}
  \Xi(x,y,z) = \infnorm{r(x,y,z)}^{\xi},
\end{equation}
where $\xi\le 1$, and $r(\cdot)$ is the KKT residual.
By default $\xi = 0.5$.

The function $\Psi(\hat{z}, x,\tau)$ is defined such that, in a local linear model of the problem, the iterates stay in the positive quadrant.
Recall that the IPM can be interpreted as solving a sequence of barrier problems with log barrier $b(c) = \mu\log(c)$.
The barrier problem for given algorithmic parameters $\mu$, $\tau$, and $\delta$, has the objective:
\begin{equation}
  f(x) - \mu\log(g(x)) -\mu\log(x_1+\delta_1) -\mu\log(x_2 + \delta_2) -\mu\log(\tau - x_1x_2).
\end{equation}
For the following analysis we will, without loss of generality, fix that $x_2 \gg 0$ such that $\mu\log(x_2 + \delta_2) \approx 0$.
To simplify the exposition we treat the case where $x = (x_1,x_2)\in \R^2$.
With these assumptions we can simplify our barrier problem objective to
\begin{equation}
  \label{eq:obj}
  f(x) - \mu\log(g(x)) - \mu\log(x_1+\delta_1) - \mu\log(\tau - x_1x_2).
\end{equation}
Note that taking the derivative w.r.t $x_1$ of this objectives yields:
\begin{equation}
  \label{eq:1dkkt}
  \nabla_{x_1}f(x) - \frac{\mu}{g(x)}\nabla_{x_1}g(x) - \frac{\mu}{x_1+\delta_1} - \frac{\mu}{\tau - x_1x_2}x_2.
\end{equation}
We observe that the last two terms in \Cref{eq:1dkkt}, $\frac{\mu}{x_1+\delta_1} - \frac{\mu}{\tau - x_1x_2}x_2$ are exactly the estimate of the MPCC multiplier $\hat{\zcc}_1$.
When close to the solution we can replace the first to terms in \Cref{eq:obj} with slope $\hat{\zcc}_1$ which yields:
\begin{equation}
  F_{\mathrm{lin}}(x_1) = \hat{\zcc}_1 x_1 - \mu\log(x_1+\delta_1) - \mu\log(\tau - x_1x_2).
\end{equation}

Given this approximation we want to develop a function $\Xi$ such that:
\begin{itemize}
\item Iterates do not leave the positive quadrant, due to the Scholtes relaxation's log barrier.
\item If $\hat{\zcc}_1 < 0$, choose a $\delta_1$ which is large enough to reduce the diagonal contribution of the lower bound log barrier term and drives $x_1$ to zero.
\end{itemize}
Therefore, it is useful to characterize the conditions on $\delta_1$ and $\tau$ which lead to a minimum at $x_1 =0$.
This holds true when:
\begin{equation}
  \nabla_{x_1} F_{\mathrm{lin}}(x_1) = \hat{\zcc}_1 - \frac{\mu}{\delta_1} - \frac{\mu x_2}{\tau} = 0.
\end{equation}
Solving this equation for $\delta_1$ yields:
\begin{equation}
  \delta_1 = \frac{\tau\mu}{\mu x_2 + \tau \hat{\zcc}_1}.
\end{equation}
These expressions only make sense when the denominator is strictly positive, and otherwise when $\hat{\zcc}_1 < 0$ the minimum is \textit{always} at $x_1 > 0$ and any lower bound relaxation $\delta_1$ can be used.
Therefore we choose:
\begin{equation}
  \label{eq:Psi_def}
  \Psi(\zcc,x,\tau) =
  \begin{cases}
    \frac{\tau\mu}{\mu x_2 + \tau \zcc} & \zcc \ge 0\ \textrm{and}\ \mu x_2 + \tau \zcc_1 > 0\\
    \delta_{\mathrm{max}} & \mathrm{else}
  \end{cases},
\end{equation}
with $\delta_{\mathrm{max}} = 10^{-4}$ as default in our implementation.

\subsection{Homotopy Parameter Updates}\label{sec:madnlpc:homotopy_update}
One key ingredient for efficiently solving nonlinear programs using interior-point methods is efficient selection of the barrier parameter $\mu$ at every iteration~\cite{Nocedal2009}.
In this section we discuss several update strategies which we implement for the now two homotopy parameters in our algorithm $\mu$ and $\tau$.
In particular, updates of the upper-level complementarity relaxation $\tau$ play an important role in the practical performance of the algorithm.

We implement several approaches to driving the relaxation parameter $\tau$ to zero.
The interior-point approach in~\cite{Raghunathan2005} ties the value of $\tau$ directly to $\mu$ by a relation $\tau(\mu):\R\rightarrow\R$ (in their case $\tau(\mu)=\mu$ is found to be optimal for the benchmarked problem set).
In the prior section on the barrier update, we made the assumption that we set $\tau^k = \mu^k$ at every iteration.
Our implementation generalizes this to $\tau^k = \alpha_{\tau}(\mu^k)^{\beta_\tau}$, which allows the implementation of e.g. $\tau_k = \sqrt{\mu^k}$ and $\tau^k = 0.1\mu^k$, and defaults to $\alpha_\tau=1$ and $\beta_\tau=1$.

However, it is not clear that this is the best choice for several reasons.
First, if a monotone $\mu$ update is used this method inherits the strong sensitivity to the initial choice of $\mu$ for a given problem.
A second effect of this choice is the stronger attraction of the algorithm to worse local minima.
We conjecture that this is caused by the low relative accuracy to which we solve the early relaxed barrier problem, taking less advantage of the convexifying effect of the relaxation.
An example of this is given in \Cref{ex:homotopy_update}.

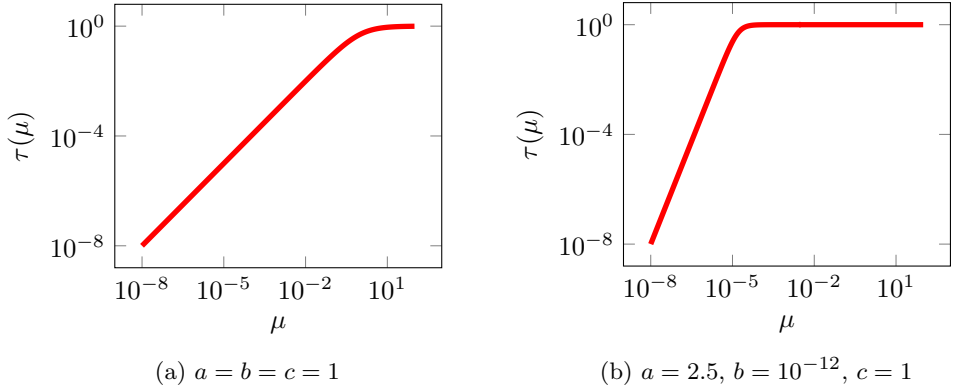
\begin{figure}[t]
  \centering
  \begin{subfigure}[b]{0.5\textwidth}
    \begin{tikzpicture}
    \begin{axis}[
      xmode=log,
      ymode=log,
      xlabel = {\(\mu\)},
      ylabel = {\(\tau(\mu)\)},
      width=.9\textwidth,
      ]
      \addplot [
      domain=1e-8:100,
      samples=1000,
      color=red,
      line width=2pt,
      ]
      {x/(x+1)};
    \end{axis}
  \end{tikzpicture}
  \caption{$a=b=c=1$}
  \end{subfigure}%
  ~
  \begin{subfigure}[b]{0.5\textwidth}
    \begin{tikzpicture}
    \begin{axis}[
      xmode=log,
      ymode=log,
      xlabel = {\(\mu\)},
      ylabel = {\(\tau(\mu)\)},
      width=.9\textwidth,
      ]
      \addplot [
      domain=1e-8:100,
      samples=1000,
      color=red,
      line width=2pt,
      ]
      {x^(2.5)/(x^(2.5)+1e-12)};
    \end{axis}
  \end{tikzpicture}
  \caption{$a=2.5$, $b=10^{-12}$, $c=1$}
\end{subfigure}
\caption{Illustration of the function~\cref{eq:rolloff_rule} for different parameter values.}
\label{fig:rolloff_plots}
\end{figure}

As such we use by default an alternative update rule, the curve for which can be found in \Cref{fig:rolloff_plots} and is defined by the expression:
\begin{equation}
  \label{eq:rolloff_rule}
  \tau(\mu)=\frac{c\mu^a}{\mu^a+b},
\end{equation}
with three parameters $a,b,c\in \R$.
The parameter $c$ defines the maximum relaxation, $b$ can be roughly interpreted as the accuracy to which the initial relaxation is solved, and $a$ as a linear slope which governs the updates after that initial phase.
Through benchmarking we observe that a good choice of parameters for this rule is $a=2.0$, $b=10^{-6}$, and $c=1$.

Additionally, we implement an alternative adaptive method for choosing $\tau^k$, inspired by the \textit{LOQO} rule~\cite{Shanno2000}:
\begin{equation*}
  \sigma_{\tau} = \gamma_{\tau}\min \paren{(1-r_\tau)\frac{1-\xi_\tau}{\xi_\tau}, 2}^3,\quad \xi_\tau = \frac{\ncc\min(X_1x_2)}{\transp{x_1}x_2},
\end{equation*}
which is used to update via:
\begin{equation*}
  \tau^{k+1} = \sigma_\tau\frac{\transp{{x_1^k}}x_2^k}{\ncc}.
\end{equation*}
This allows for nonmonotonicity in the behavior of the  complementarity relaxation parameter and reduces the effect of the initial choices of $\tau^0$ as well as $\alpha_\tau$ and $\beta_\tau$.

We now turn to describing the choice of update rule for the barrier parameter $\mu$.
The baseline approach for updating the barrier parameter is the static approach of Fiacco-McCormick which uses a fixed accelerating profile to drive $\mu$ to zero:
\begin{equation*}
  \mu^{k+1} = \alpha_{\mu}(\mu^{k})^{\beta_\mu},
\end{equation*}
with $\beta > 1$ and $\alpha \in (0, 1)$.
This update is applied when the barrier problem has been solved to a tolerance factor of $\mu^{k}$.
The default values taken in \madnlp\ are $\alpha_\mu = 0.2$ and $\beta = 1.5$ with $\mu^0 = 0.1$ which are shared in {\solver}.
This method has been shown to have global convergence properties for NLPs but is sensitive to initial barrier parameter $\mu^0$~\cite{Waechter2006a,Waechter2006}.

Other approaches can be used that include the so called \textit{LOQO rule} used in the \textit{eponymous} software package~\cite{Shanno2000}.
We also investigate an adaptive barrier rule~\cite{Nocedal2009} which attempts to choose a step in $\mu$ which minimizes a ``quality function'' which is generally an approximation of the norm of the KKT error in \Cref{eq:term_crit}.

Both adaptive rules generally follow the pattern of selecting a $\sigma^k > 0$ at each iteration which is then used to update the barrier parameter via
the average complementarity at the current iterate:
\begin{equation}
  \label{eq:mu_update_rule}
  \mu^{k+1} = \sigma^k\frac{ (x^{k})^\top z^{k} }{n}.
\end{equation}
Recall that $n$ here is the dimension of all primal variables $x = (x_0,x_1,x_2)$.
We extend both adaptive rules to our MPCC setting, where both $\tau$ and $\mu$ are driven simultaneously to zero.
The extension includes both the upper-level complementarity constraints ($ 0 \leq x_1 \perp x_2 \geq 0$, in the MPCC constraints), and the lower-level complementarity conditions in the KKT conditions of the relaxed MPCC:
\begin{equation} \label{eq:mpcc_mu_update_rule}
  \mu^{k} = \sigma^k\frac{({x^{k}})^\top z^k + {({x_1^k})^\top}x_2^k}{n+\ncc}.
\end{equation}
This extension clearly translates when $\tau^{k} = \mu^{k}$ is chosen as the relationship between the two complementarity relaxations.
For all of the adaptive methods described below we further use the globalization scheme implemented in {\madnlp} and described in~\cite{Nocedal2006}.
This approach can be summarized as follows: we check for sufficient progress in the KKT residual, and if none is achieved, fall back to the monotone rule.
Note, this globalization scheme is parallel to the filter line-search scheme implemented to ensure global convergence of the interior-point method.

\paragraph{A Modified LOQO rule.}
The \textit{LOQO rule}~\cite{Vanderbei1999} is defined as follows:
\begin{equation}
  \label{eq:loqo_rule}
  \sigma^k = \gamma\min \paren{(1-r)\frac{1-\xi}{\xi}, 2}^3,\quad \xi = \frac{n\min(X^kz^k)}{\transp{(x^k)}z^k},
\end{equation}
with \textit{step length} parameter $r\in (0,1)$, and scaling parameter $\gamma > 0$.
In our case, we propose the modification to this rule that incorporates the upper-level complementarity values:
\begin{equation}
  \label{eq:mpcc_loqo_rule}
  \sigma^k = \gamma\min \paren{(1-r)\frac{1-\xi}{\xi}, 2}^3,\quad \xi = \frac{(n_x+n_{cc})\min((X^kz^k,X^k_1x^k_2))}{\transp{(x^k)}z^k + \transp{(x_1^k)}x_2^k}.
\end{equation}
In {\solver}, we implement both approaches to selecting $\sigma^k$.
In our benchmarking the joint rule in \cref{eq:mpcc_loqo_rule} performed better in practice and is used by default when the user chooses the \textit{LOQO rule} option.

\paragraph{Quality-function based update.}
The second adaptive barrier method we implement in {\solver}, is inspired by the work done by Nocedal et al.~\cite{Nocedal2009}, which attempts to choose an \textit{optimal} $\sigma$, in terms of reduction in KKT error.
For simplicity of exposition, let us take $\mu(\sigma)$ to be calculated as in \Cref{eq:mpcc_mu_update_rule}, that is, including the upper-level complementarities.
We define the steps resulting from the solution of \Cref{eq:aug_kkt} for a given $\sigma$ as $\Delta x(\sigma)$, $\Delta s(\sigma)$, $\Delta y(\sigma)$, and $\Delta z(\sigma)$.
The maximum fraction-to-boundary step (i.e. satisfying \cref{eq:frac_to_boundary}) in the primal and dual variables also depend on $\sigma$ and as such we write these as $\aprmax(\sigma)$ and $\adumax(\sigma)$.
This then yields:
\begin{align*}
  x(\sigma) &= x + \aprmax(\sigma)\Delta x(\sigma),\ s(\sigma) = s + \aprmax(\sigma)\Delta s(\sigma),\\
  y(\sigma) &= y + \adumax(\sigma)\Delta y(\sigma),\ z(\sigma) = z + \adumax(\sigma)\Delta z(\sigma).
\end{align*}

As in~\cite{Nocedal2009} we choose to minimize the norm of the KKT error:
\begin{equation*}
  \Psi(x,s,y,z) = \twonorm{\nabla\Lag(x,s,y,z)}^2 + \twonorm{c(x)}^2 + \twonorm{Zx}^2 + \twonorm{X_1x_2}^2,
\end{equation*}
with the inclusion of the multiplicative residual of the complementarity variables.
This optimal reduction is calculated through the minimization of a so called ``quality-function'', which is chosen to be a proxy for the KKT error.
In principle the quality function can be taken as a design choice, however, it is important to note that in order to maintain the effectiveness of the method, additional evaluations of the Hessian of the Lagrangian, constraint Jacobian, and constraint gradient should be avoided, as these can be quite expensive.
To that end we similarly choose a linear approximation to the true KKT error:
\begin{equation}
  \label{eq:lin_q_fun}
  \begin{aligned}
    q_L(\sigma) = &(1-\adumax(\sigma))^2\twonorm{\nabla\Lag(x^k,s^k,y^k,z^k)}^2 + (1-\aprmax(\sigma))^2\twonorm{c(x^k)}^2 +\\
                  & \twonorm{(Z^k+\adumax(\sigma)\Delta Z(\sigma))(x^k+\aprmax(\sigma)\Delta x(\sigma))}^2 + \\
                  & \twonorm{(X_1^k + \aprmax(\sigma)\Delta X_1(\sigma))(x_2^k + \aprmax(\sigma)\Delta x_2(\sigma)}^2.
  \end{aligned}
\end{equation}
This linear approximation differs from \cite[Equation 4.2]{Nocedal2009}, applied to the Scholtes-relaxed NLP directly as it does not use the linearization of the Scholtes relaxation.
On the contrary, we use the norm of the exact calculation of the upper-level complementarity residual for two reasons.
(i) It is cheap to evaluate exactly, requiring only vector-vector operations.
(ii) It is more accurate to the ``goals'' of the solver, that is to minimize the residual, including that of the upper level complementarities $X_1x_2$.
As such we only use the linearization of the general nonlinear constraints and the gradient of the Lagrangian, tailoring the quality function for {\MPCCs}.

As a choice of $\mu$ and $\tau$ affects only the right hand side of the Newton step calculated by \Cref{eq:aug_kkt}, the evaluation of \Cref{eq:lin_q_fun} requires one numeric factorization, two backsolves, and some vector arithmetic.
In particular we first solve the KKT system in \Cref{eq:aug_kkt} with $\mu = 0$ and $\tau = 0$, which is related to the affine step in the classic Mehrotra predictor-corrector scheme.
This gives us the steps $\Delta x^{\mathrm{aff}}$, $\Delta s^{\mathrm{aff}}$, $\Delta y^{\mathrm{aff}}$, and $\Delta z^{\mathrm{aff}}$.
Then we calculate the centering step, $\Delta x^{\mathrm{cen}}$, $\Delta s^{\mathrm{cen}}$, $\Delta y^{\mathrm{cen}}$, and $\Delta z^{\mathrm{cen}}$, by solving \Cref{eq:aug_kkt} with the right hand side corresponding to:
\begin{alignat*}{2}
  r_i &= 0,\ &i=1,\ldots 5,\\
  r_i &=\frac{\transp{{x^{k}}}z^k + \transp{{x_1^k}}x_2^k}{n+\ncc},\ &i=6\ldots 10.
\end{alignat*}
Note that by linearity we can calculate $\Delta x(\sigma) = \Delta x^{\mathrm{aff}} + \sigma\Delta x^{\mathrm{cen}}$.
The same can be done for steps $\Delta s(\sigma)$, $\Delta y(\sigma)$, and $\Delta z(\sigma)$.
We then approximately minimize \Cref{eq:lin_q_fun} using the golden-search approach described in~\cite{Nocedal2009}.

\subsection{Some illustrative examples}
\begin{figure}[t]
  \centering
  \begin{subfigure}{0.45\linewidth}
    \centering
    \includegraphics[width=\textwidth]{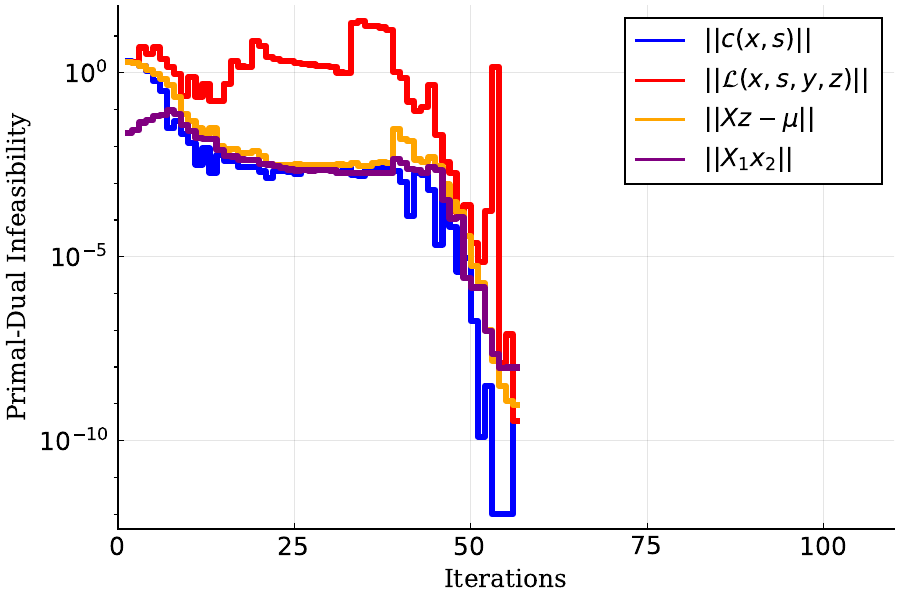}
    \caption{Feasibility, proportional rule.}
  \end{subfigure}
  \begin{subfigure}{0.45\linewidth}
    \includegraphics[width=\textwidth]{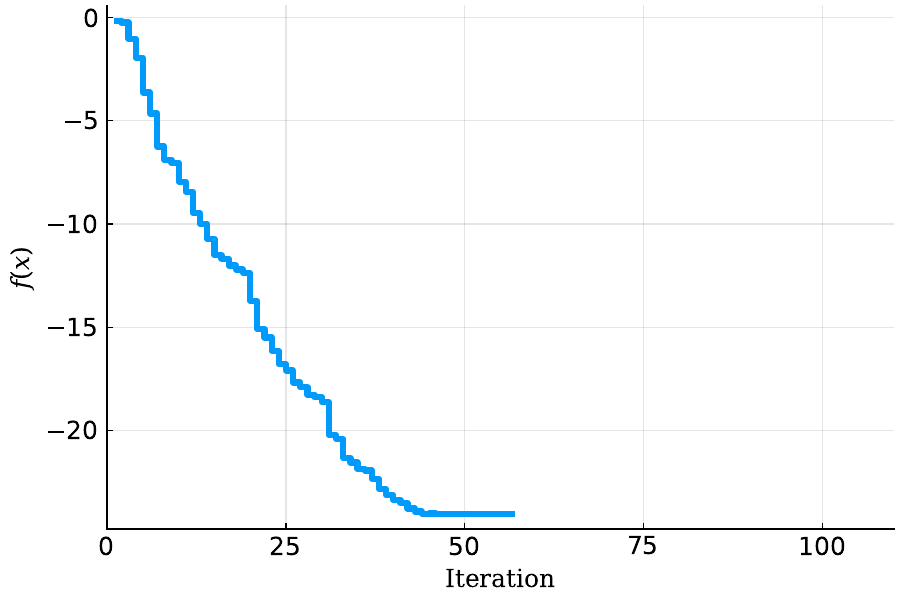}
    \caption{Objective, proportional rule.}
  \end{subfigure}
  \\
  \begin{subfigure}{0.45\linewidth}
    \centering
    \includegraphics[width=\textwidth]{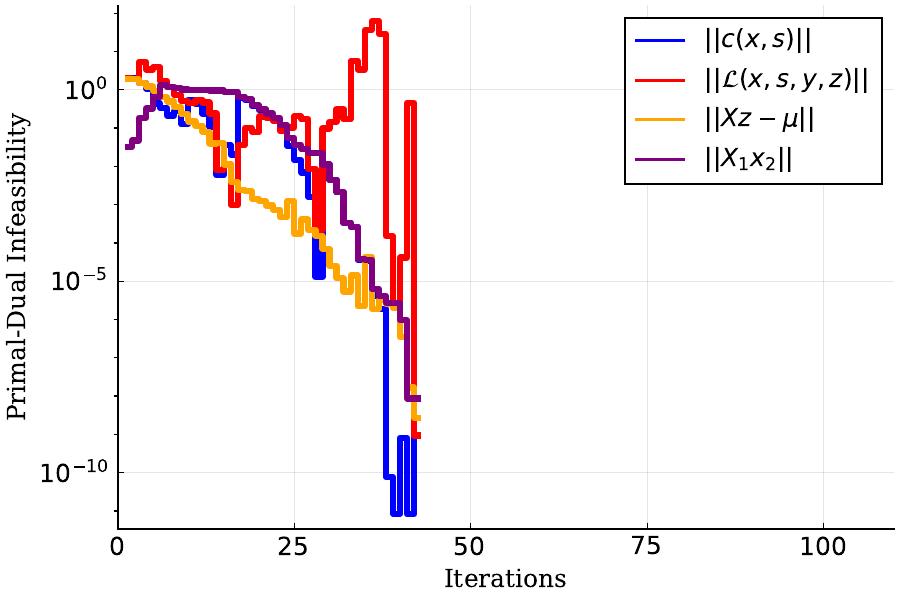}
    \caption{Feasibility, alternative rule.}
  \end{subfigure}
  \begin{subfigure}{0.45\linewidth}
    \includegraphics[width=\textwidth]{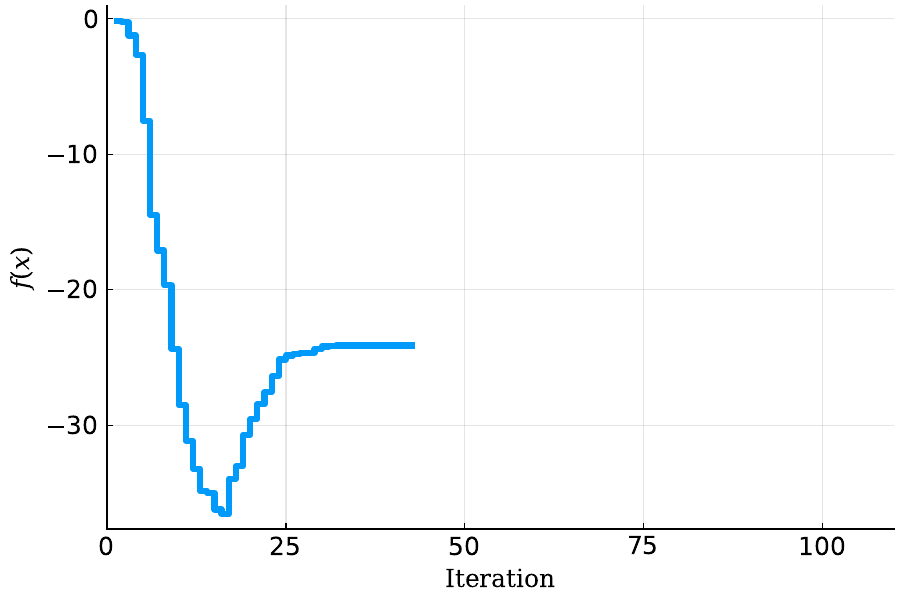}
    \caption{Objective, alternative rule.}
  \end{subfigure}
  \\
  \begin{subfigure}{0.45\linewidth}
    \centering
    \includegraphics[width=\textwidth]{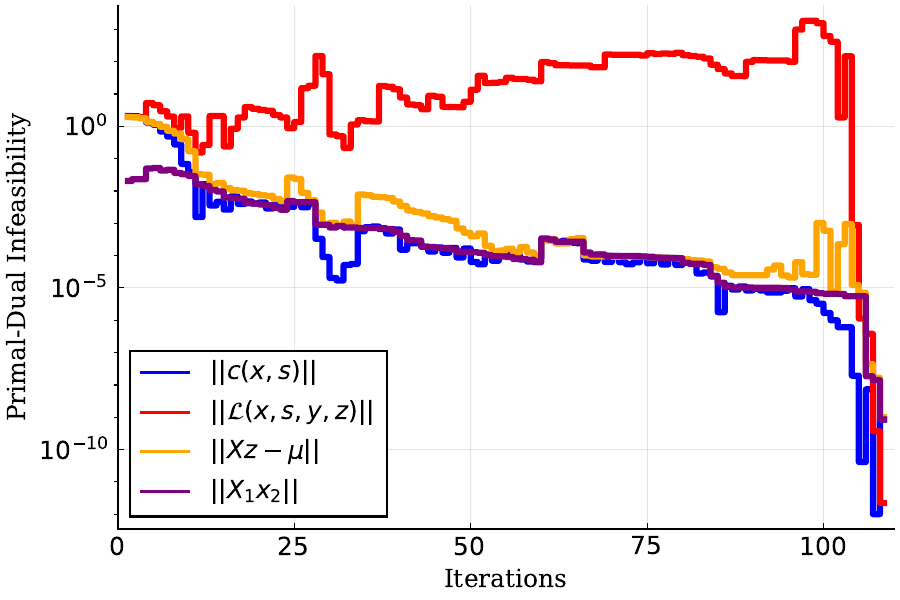}
    \caption{Feasibility, quality function update rule.}
  \end{subfigure}
  \begin{subfigure}{0.45\linewidth}
    \includegraphics[width=\textwidth]{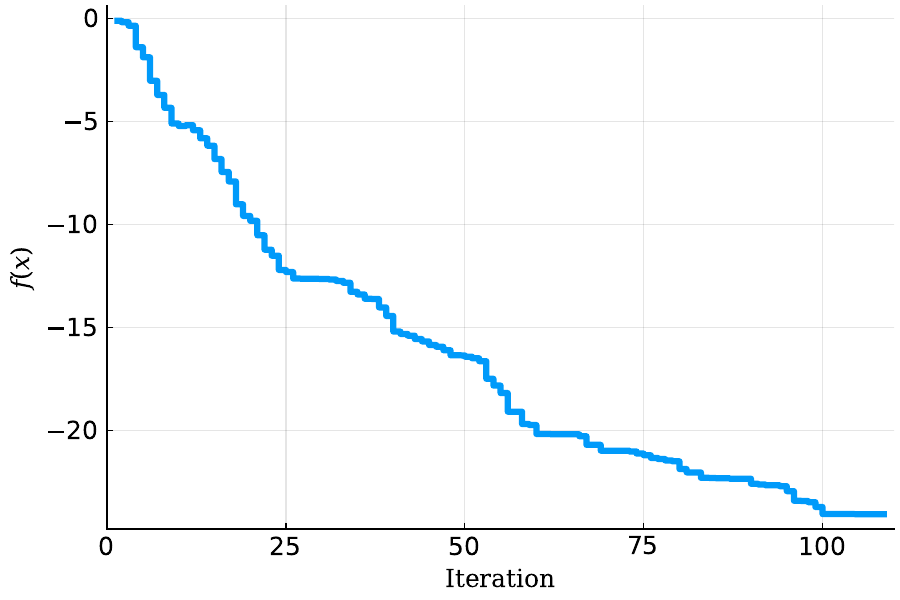}
    \caption{Objective, quality function update rule.}
  \end{subfigure}
  \caption{Solver behavior on \texttt{qpec-200-2} using different parameter update rules.}
  \label{fig:qpec_200}
\end{figure}

\begin{figure}[t]
  \centering
  \begin{subfigure}{0.45\linewidth}
    \centering
    \includegraphics[width=\textwidth]{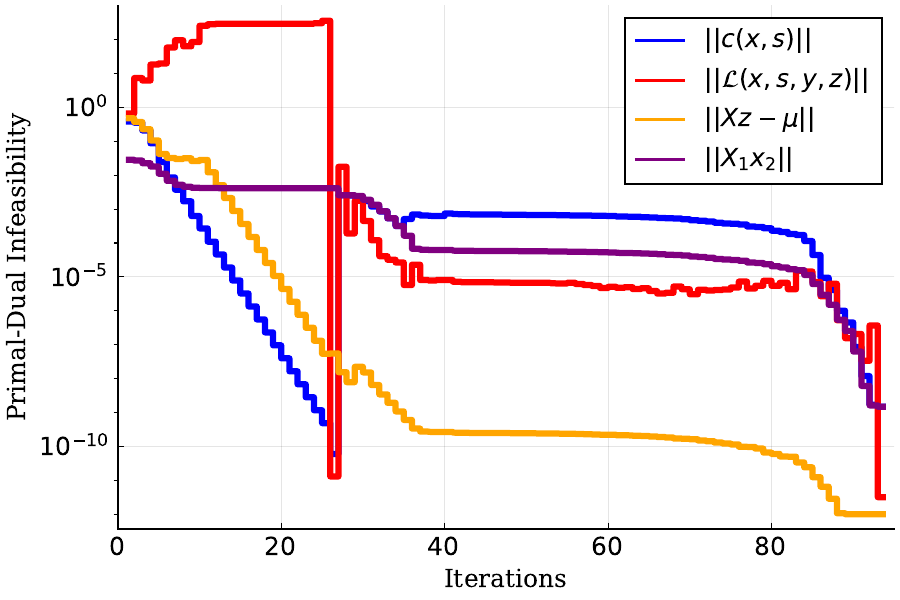}
    \caption{Feasibility, proportional rule.}
  \end{subfigure}
  \begin{subfigure}{0.45\linewidth}
    \includegraphics[width=\textwidth]{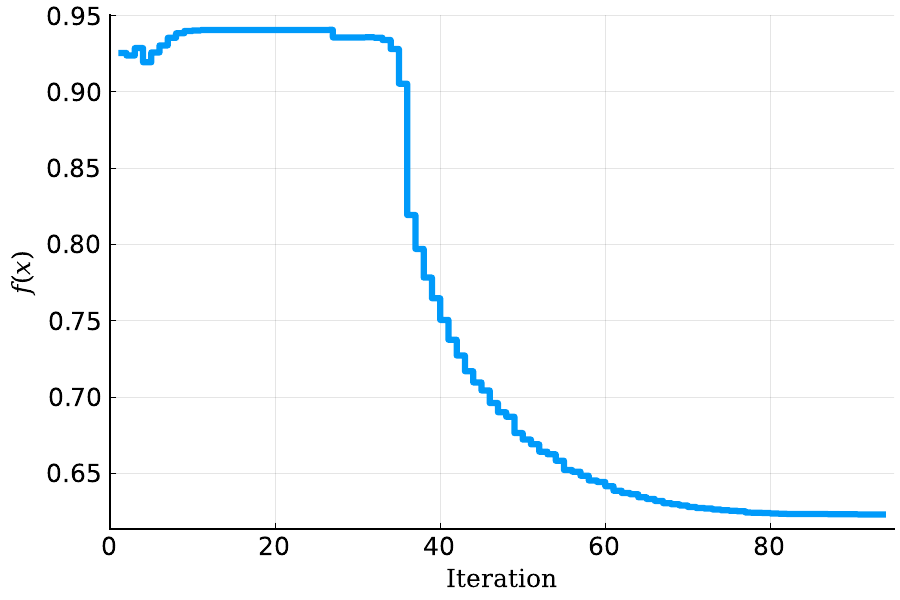}
    \caption{Objective, proportional rule.}
  \end{subfigure}
  \\
  \begin{subfigure}{0.45\linewidth}
    \centering
    \includegraphics[width=\textwidth]{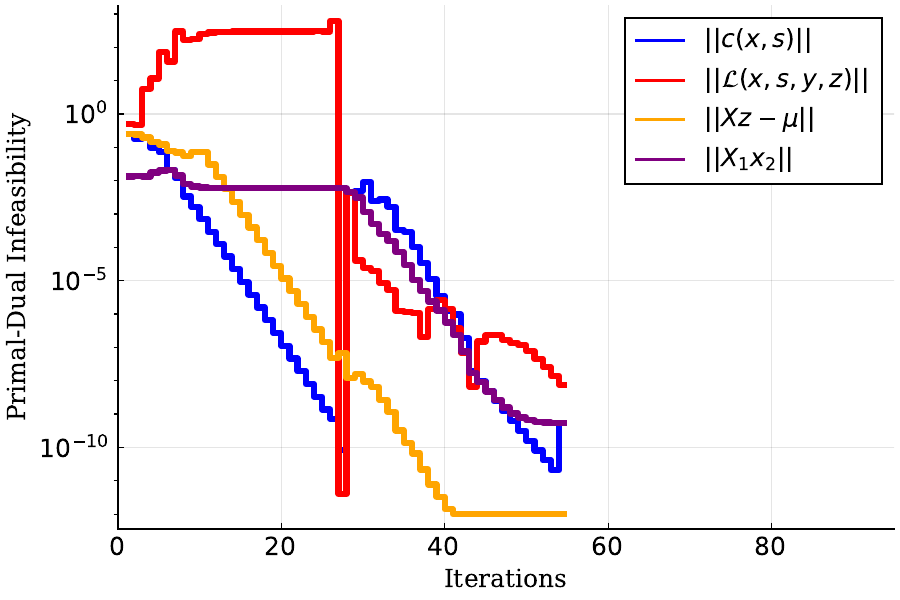}
    \caption{Feasibility, quality function update rule.}
  \end{subfigure}
  \begin{subfigure}{0.45\linewidth}
    \includegraphics[width=\textwidth]{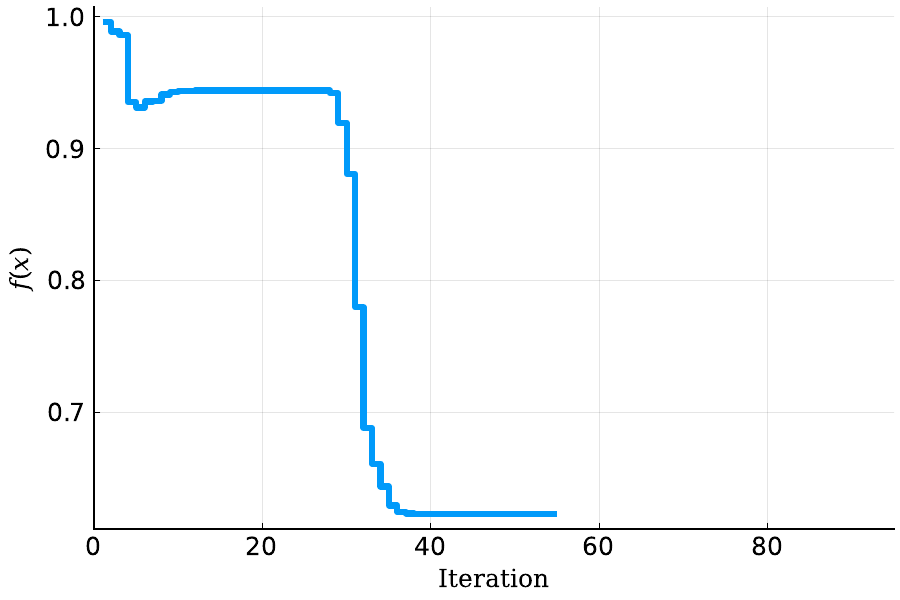}
    \caption{Objective, quality function update rule.}
  \end{subfigure}
  \caption{Solver behavior on \texttt{pack-comp1c-16} using different update rules.}
  \label{fig:pack_comp}
\end{figure}

Finally, we provide two illustrative examples of the algorithmic choices proposed above which are representative of their effects.
We first compare the homotopy parameter updates and describe their beneficial and detrimental effects on solver performance and conclude with an example of the efficacy of the proposed ``Q-regularization'' scheme, particularly in the case of quadratic programs with complementarity constraints.
\begin{example}[A comparison of homotopy update strategies]\label{ex:homotopy_update}
  Here we compare solver behavior on several problems with the proportional rule $\tau = \mu$, the alternative rule, as well as the quality function based update rule.
  The first problem we regard is \texttt{qpec-200-2}~\cite{Jiang1999} from the {\MacMPEC}~\cite{Leyffer2000} problem set.
  This problem has many local minima of different quality: we demonstrate that the alternative update rule's benefit of finding better local minima.
  We provide plots of the objective and infeasibility norms for both the proportional and the alternative rule in \Cref{fig:qpec_200}.
  We note that these two methods find \textit{different} local minima with different complementarity active sets, with the minimum found by the alternative update rule having a lower objective value ($-2.411\times 10^{1}$ versus $-2.407\times 10^{1}$).
  This behavior is quite valuable, particularly in optimal control problems which may exhibit many qualitatively bad B-stationary local minima~\cite{Nurkanovic2020}.
  We also plot the solver trajectory using the joint updating quality function based approach in \Cref{fig:qpec_200}.
  The minimum found by this method once again has a different complementarity active set and a worse objective ($-2.406\times 10^{1}$) and takes many more iterations to converge.

  To show the usefulness of the quality function update rule on problems with fewer local minima: we now use the problem \texttt{pack-comp1c-16}, which comes from~\cite{Outrata2013}.
  In \Cref{fig:pack_comp}, the reader can find a comparison of the behavior of the solver using the proportional update rule and the joint quality function update rule, respectively.
  In this case, both methods find the same local minimum, but the quality function based approach takes better steps in the homotopy parameters, which leads to a nearly factor of two reduction in the number of iterations it takes to converge.
  We observe this behavior in many cases where progress using the joint quality function based approach prevents long periods of stalled improvement.

  These examples effectively show that the choice of update rules for $\mu$ and $\tau$ is non-trivial, and it is difficult to say with certainty that one is clearly the best option for all problems.
  On the other hand, understanding something of the problem's structure can lead to a more informed choice of options.
  In this article, we seek to also provide some guidance on what options to choose when using {\solver}.
\end{example}

\begin{figure}[t]
  \centering
  \includegraphics[width=\textwidth]{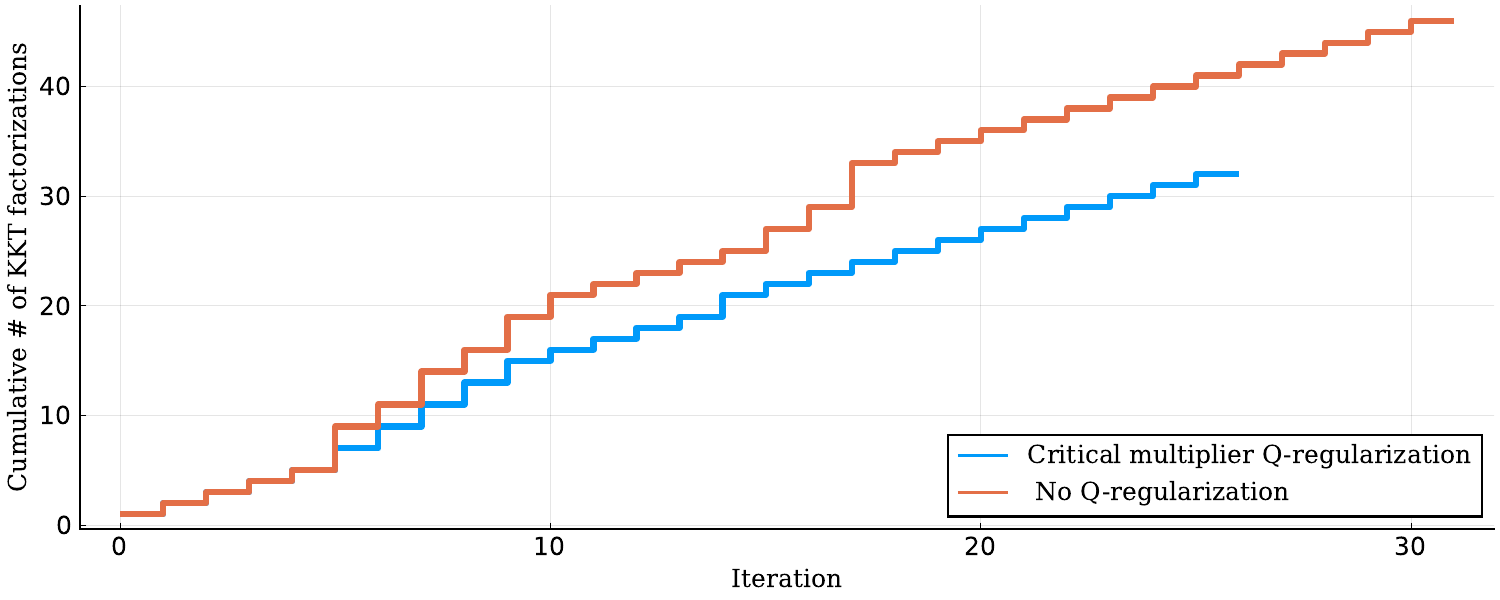}
  \caption{Cumulative KKT Matrix factorizations at each iteration, comparing critical multiplier regularization with no regularization.}
  \label{fig:qpec_100_reg}
\end{figure}

\begin{example}[Q-regularization]
  To illustrate the usefulness of the Q-regularization scheme described in \Cref{sec:madnlpc:hessian_regularization} we compare the behavior of the {\solver} relaxation algorithm with and without the critical multiplier regularization scheme.
  This technique is particularly useful for QPCC problems where the original objective is a convex quadratic, such as \texttt{qpec-100-2} from the {\MacMPEC} test set.
  In \Cref{fig:qpec_100_reg} we show the number of KKT matrix factorizations over the number of iterations.
  Note, that when using the critical multiplier scheme we take at most two KKT factorizations per iteration, whereas when using the inertia correction scheme, up to four factorizations are needed to identify a regularization sufficiently large.
  This example also exhibits an improvement in the number of iterations to convergence when using the critical multiplier scheme.
  While we frequently observe this to be the case in practice, there is in principle no guarantee that the Q-regularization schemes improve iteration counts over the classic inertia correction approach.
\end{example}

\subsection{On convergence properties}
The convergence properties of relaxation homotopy methods for {\MPECs} have been studied extensively in the literature and as such we will not reprise this theory in its entirety.
The local convergence properties for this Scholtes relaxation homotopy can be found in~\cite{Raghunathan2005}, which provides a \textit{q-quadratic} convergence for points sufficiently close to S-stationary points in conjunction to the Vicente-Wright regularization scheme, under the strict complementarity assumption in~\cite{Scheel2000} and MPCC-MFCQ.
The convergence of approximate solutions of the Scholtes relaxation was studied in~\cite{Kanzow2015}, which suggests that unlike many other relaxations $\epsilon$-stationary points of the relaxed NLP converge to M-Stationary points of the original problem, under MPCC-MFCQ.
Our modifications to the selection of the homotopy parameters do not impact the local convergence properties of the algorithm: close to the solution the parameters $\mu$ and $\tau$ behavior is consistent with the theory in~\cite{Raghunathan2005}.
Global convergence of our algorithm follows trivially from the global convergence properties of filter line search~\cite{Waechter2005} and the globalization arguments from~\cite{Raghunathan2005,Nocedal2009}.

\subsection{On initialization}\label{sec:initialization}
Initialization complexity is a long standing pain-point for interior-point methods, and the algorithms we describe are no different.
We observe that in practice (unless a good guess for the complementarity variables is known a priori) it can help to initialize the complementarities in a ``centered'' fashion.
As such we implement the option for the relaxation based algorithm to initialize the complementarity variables and the corresponding Scholtes slack ($x_1$, $x_2$ and $s$ in \Cref{eq:madnlpc_nlp}) automatically to a feasible point.
This is done for a fixed $\tau^0$ and a $k_{\mathrm{cen}}\in (0,1)$, which is set to $0.5$ by default:
\begin{align*}
  x_1 &= x_2 = \sqrt{k_{\mathrm{cen}}\tau^0},\\
  s &= \sqrt{2(1-k_{\mathrm{cen}}\tau^0)}.
\end{align*}
For many problems this relaxed complementarity feasible initialization often, though not always, improves performance over standard bound pushing initializations which are described in~\cite[Section 3.6]{Waechter2006}.

\section{An interior-penalty algorithm for MPCCs}
\label{sec:penalty}
In this section we describe the implemented interior penalty method in the style of Leyffer et al.~\cite{Leyffer2006}, where a sequence of penalty problems are solved with the penalty parameter updated dynamically either at each iteration or at the solution of each barrier problem.
We focus on describing a novel Hessian regularization, as well as several dynamic penalty algorithms.
The high level overview of the algorithms can be found in \Cref{alg:interior_penalty}.

The methods described in this section solve the following ``exact-penalty'' NLP derived from \Cref{eq:mpec}:
\begin{mini!}[4]
  {\substack{x\in\R^{n}}}
  {f(x) + \rho \transp{x_1}x_2}
  {\protect\label{eq:interior_penalty_nlp}}
  {}
  \addConstraint{c(x)}{= 0 \; ,}
  \addConstraint{x}{\ge 0 \; .}
\end{mini!}

\begin{algorithm}[b]
  \caption{Interior penalty algorithm for {\MPCCs}.}
  \label{alg:interior_penalty}
  \begin{algorithmic}[1]
    \State $k = 0$
    \State $\rho = \rho_0$
    \While{\eqref{eq:term_crit} is not satisfied}
      \State $\mu^k = \mathrm{update\_mu}(x^k,s^k,y^k,z^k)$
      \State $\rho = \mathrm{update\_rho}(x^k,s^k,y^k,z^k,\mu^k)$
      \State Evaluate KKT matrix $K^k$ and KKT residual $r^k$.
      \State Calculate search direction $d^k$ by solving $K^kd^k = b^k$.
      \State Calculate step lengths $\apr$ and $\adu$ using a filter line-search scheme.
      \If{step sizes are too small or search direction calculation fails}
        \State Begin restoration phase for fixed $\rho$ a la \cite{Waechter2006}.
      \Else
        \State $x^{k+1} = x^k + \apr d_x^k$, $s^{k+1} = s^k + \apr d_s^k$
        \State $y^{k+1} = y^k + \adu d_y^k$, $z^{k+1} = z^k + \adu d_z^k$
        \State $k=k+1$
      \EndIf
    \EndWhile
  \end{algorithmic}
  \vspace{-0.1cm}
\end{algorithm}
\subsection{Termination criteria}\label{sec:ell1:term_crit}
In this method we use identical termination criteria as in \Cref{eq:term_crit}.
In practice we observe that penalty methods can often achieve much better complementarity residuals with comparatively less numerical difficulties, when problems are easy.
However, to make equivalent comparisons viable between our two algorithms we use the same default $\epstol=10^{-8}$ for these algorithms.

\subsection{Step calculation}\label{sec:ell1:step_calculation}
The (already relaxed) KKT conditions for \Cref{eq:interior_penalty_nlp} can be written as:
\begin{alignat*}{2}
  \nabla_{x_0}f(x) - \nabla_{x_0}h(x)y_h \eqqcolon r_1 &= 0,\\
  \nabla_{x_1}f(x) + \rho x_2 - \nabla_{x_1}h(x)y_h + X_2y_s - z_1 \eqqcolon r_2 &= 0,\\
  \nabla_{x_2}f(x) + \rho x_1 - \nabla_{x_2}h(x)y_h + X_1y_s - z_2 \eqqcolon r_3 &= 0,\\
  c(x) \eqqcolon r_4 &= 0,\\
  X_0z_0 - \mu e \eqqcolon r_5 &= 0,\\
  X_1z_1 - \mu e \eqqcolon r_6 &= 0,\\
  X_2z_2 - \mu e \eqqcolon r_7 &= 0,\\
  x,z&> 0.
\end{alignat*}
The corresponding augmented KKT system solved at every iteration of the exact penalty method is then written as:
\begin{equation}
  \label{eq:aug_kkt_pen}
  Kd = -r
\end{equation}
with
\begin{equation*}
  K =
  \begin{bmatrix}
    W + Q_\rho &\transp{B}\\
    B & 0
  \end{bmatrix},\quad
  d =
  \begin{bmatrix}
    \Delta x\\
    \Delta y
  \end{bmatrix}, \quad
  r=
  \begin{bmatrix}
     r_{x}\\
     r_5
  \end{bmatrix},
\end{equation*}
\begin{equation*}
  W = \nabla^2_x(f(x)+\transp{c(x)}y_c)=
  \begin{bmatrix}
    W_{00} & W_{01} & W_{02}\\
    W_{00} & W_{11} & W_{12}\\
    W_{20} & W_{21} & W_{22}
  \end{bmatrix},\quad
  Q_\rho =
  \begin{bmatrix}
    \Sigma_0 & 0 & 0\\
    0 & \Sigma_1 & \rho I\\
    0 & \rho I & \Sigma_2
  \end{bmatrix},
\end{equation*}
\begin{equation*}
  B =
  \begin{bmatrix}
    J_0& J_1& J_2
  \end{bmatrix},\quad
  r_x =
  \begin{bmatrix}
    r_1+\inv{X}_0r_7\\
    r_2+\inv{X}_1r_8\\
    r_3+\inv{X}_2r_9
  \end{bmatrix},\quad
\end{equation*}
Note that this KKT matrix has a nearly identical shape to that solved in the relaxation based algorithm, without the final row and column coming from the Scholtes constraints, and the slacks.
The multiplier for the absent constraints is replaced by the algorithmic penalty parameter $\rho$.

\subsection{Hessian regularization}\label{sec:ell1:hess_regularization}
We implement a similar ``$Q$-regularization'' technique as described in \Cref{sec:madnlpc:hessian_regularization}.
Note that the structure of $Q_\rho$ is identical to that of $Q$ in that section and therefore both the regularization techniques described there can be used for the penalty algorithm.

\subsection{Homotopy parameter update}\label{sec:ell1:parameter_update}
We implement the two different methods of updating the penalty parameter described in~\cite{Leyffer2006}.
The basic update rule we use in both cases is $\rho^k = \alpha_\rho \rho^k$ with $\alpha_\rho > 1$.
The first uses a static update rule which applies a fixed update rule when the current barrier problem is solved to sufficient accuracy.
This method may, however, fail in some cases when the penalty problem is unbounded for a given penalty parameter $\rho$.

The second ``dynamic'' approach tracks the current complementarity residual and if there is insufficient decrease, or even an increase, in the residual $\rho$ is increased dynamically, to avoid the unbounded behavior that the static method may have.
Updates are triggered when
\begin{align*}
  c(x^k)&\ge  (\mu^k)^\gamma,\\
  \transp{{x_1^k}}x_2^k &\ge \eta_{\mathrm{pen}}\frac{\sum_{i=k-h}^k\transp{{x_1^i}}x_2^i}{h}.
\end{align*}
with $h$ being the length of filter used to calculate the required decrease in complementarity, and $\eta_{\mathrm{pen}}$ the required decrease.
The basic update rule we use in both cases is $\rho^k = \alpha_\rho \rho^k$ with $\alpha_\rho > 1$.

When using the ``dynamic'' approach we also support all of the $\mu$ update methods described in \Cref{sec:madnlpc:homotopy_update}.
In contrast to the Scholtes relaxation based approach it is nontrivial to update the penalty $\rho$ in a similar way to $\tau$, as it appear on the left hand side of the KKT system.
This means that calculating an ``optimal'' update to $\rho$, would require numerically refactorizing the KKT matrix for each probed point, which would be computationally impractical.

\subsection{On convergence properties}
The global convergence properties of penalty methods have also been studied extensively~\cite{Anitescu2007,Leyffer2006}.
In particular as our algorithm corresponds to the algorithms presented in~\cite{Leyffer2006}, the same convergence properties apply, namely that \Cref{alg:interior_penalty}, converges to C-Stationary points of the {\MPCC} under the assumption of {\MPCC-LICQ}.


\section{Crossover to complementarity active set}
\label{sec:crossover}

Recently, Nurkanovi\'c and Leyffer have proposed an active-set-based method~\cite{Nurkanovic2025}, which uses Linear Programs with Complementarity Constraints (\LPCCs) to project to feasible points and verify B-Stationarity.
We provide a high-level overview of the proposed algorithm in \Cref{alg:mpecopt}.
Here, we discuss several aspects of the implementation of this algorithm in {\solver}.
First, we discuss the complementarity active-set method which provides the basis of the crossover.
Second, we discuss the several procedures for the initial complementarity active-set identification for nearly feasible points.
Finally, we conclude with a discussion of the implemented methods for solving the {\LPCC} in line 5 of \Cref{alg:mpecopt}.
Since the $\BNLP$ in line 6 of \Cref{alg:mpecopt} is a standard NLP, it is solved using the filter line-search method implemented in \madnlp~\cite{Shin2024}, on top of which {\solver} is built. 

\subsection{{\MPCC} active-set method}
This method can be divided into two phases: a Phase I finding a feasible $\BNLP$, and then a Phase II which solves a finite sequence of $\BNLP$s until a B-stationary point is found.
As there are a myriad of methods to find an initial feasible branch for the Phase I of this algorithm, details of which can be found in~\cite{Nurkanovic2025}, we focus on describing the Phase II.

The Phase II of this method can be shown to globally converge to a B-Stationary point once a feasible branch has been found.
In particular, we need a partition of the complementarity constraints into:
\begin{align*}
  \I_1 &\coloneqq \Set*{i}{x_{1,i} q\ge 0, x_{2,i} = 0}\\
  \I_2 &\coloneqq \Set*{i}{x_{1,i} = 0, x_{2,i} \geq 0},
\end{align*}
which result in a feasible $\BNLP$.
We provide a high level pseudocode of this algorithm in \Cref{alg:mpecopt}.
The $\textrm{LPEC}(x,\Delta)$ solved at each iteration of the algorithm is the so-called ``trust-region'' \LPCC:
\begin{mini!}[4]
  {\substack{d\in \R^n}}
  {\nabla f(x)^\top d}
  {\label{eq:tr_lpec}}
  {}
  \addConstraint{0}{= c(x)+ \nabla c(x)^\top d}{}
  \addConstraint{0}{=x_{1,i} + d_{1,i},}{\qquad\qquad\forall  i \in \mathcal{I}_{0+}(x)}
  \addConstraint{0}{=x_{2,i} + d_{2,i},}{\qquad\qquad\forall  i \in \mathcal{I}_{+0}(x)}
  \addConstraint{0}{\le x_{1,i} + d_{1,i}\perp x_{2,i} + d_{2,i}\ge 0,}{\qquad\qquad\forall i \in \mathcal{I}_{00}(x)}
  \addConstraint{0}{\le \Delta - ||d||_{\infty},\label{eq:tr_lpec:tr}}
\end{mini!}
which is equivalent to \Cref{eq:lpec_reduced_theory}, with the addition of an infinity-norm trust-region on the step size (\ref{eq:tr_lpec:tr}).
The trust-region ensure the \LPCC\ is bounded an non-optimal points, but it does not change the solution at B-stationary points, cf.~\cite{Nurkanovic2025}.

The second subproblem solve in phase II is denoted by $\mathrm{BNLP}{(\I_1,\I_2)}$. 
For a given partition of the complementarities $\I_1,\ \I_2$ is defined as follows:
\begin{mini!}[4]
  {\substack{x\in \R^n}}
  {f(x)}
  {\label{eq:bnlp}}
  {}
  \addConstraint{0}{=c(x)}
  \addConstraint{0}{=x_{2,i},\ 0\le x_{1,i},\; }{\forall i \in \I_1}
  \addConstraint{0}{=x_{1,i},\ 0\le x_{2,i},\; }{\forall i \in \I_2.}
\end{mini!}
\begin{algorithm}
  \caption{Complementarity active-set method for {\MPECs}}
  \label{alg:mpecopt}
  \begin{algorithmic}[1]
    \Input
    \Desc{$\mathcal{I}^0_1,\mathcal{I}^0_2\in\mathbb{B}^{\ncc}$}{A feasible branch decomposition for \MPEC~\eqref{eq:mpec}}.
    \Desc{$x^0$}{A feasible point for for $\mathrm{BNLP}{(\I^0_1,\I^0_2)}$}.
    \Desc{$\Delta^{0,0}$}{An initial trust region radius.}
    \State $k = 0$
    \While{$d^{k,l} \neq 0$} \Comment{\textcolor{gray}{major iterations}}
    \State $x^{k,0} = x^k$
    \For{$l=0,\ldots$} \Comment{\textcolor{gray}{minor iterations}}
    \State Solve $\LPEC(x^k,\Delta^{k,l})$ to obtain $d^{k,l}$
    \State Solve $\BNLP(\I_1(x^k+d^{k,l}),\I_2(x^k+d^{k,l}))$~\eqref{eq:bnlp} to obtain $x^{k,l}$
    \If{ $f(x^{k,l}) < f(x^{k})$ }
    \State Set $x^{k+1} = x^{k,l}$, increase trust region radius $\Delta^{k,l+1}$, and \textbf{break}
    \Else
    \State Reduce trust-region radius $\Delta^{k,l+1}$
    \EndIf
    \EndFor
    \EndWhile
  \end{algorithmic}
  \vspace{-0.1cm}
\end{algorithm}

\subsection{Complementarity active-set identification}
In practice, one may not have immediately an {\MPEC}-feasible point available, as is required in \Cref{alg:mpecopt}.
Thus, it is necessary to identify an initial feasible complementarity active set from a (possibly  mildly) infeasible initial point $x^0$.
As shown in~\cite[Theorem 4.1]{Nurkanovic2025}, one simply needs the point $x^0$ to be sufficiently close to feasible set of the {\MPEC}, then a feasible point of the \LPCC\ constructed at this point to identify a feasible branch NLP.
As such, an active-set identification scheme we implement involves solving the following trust region problem, $\PROJLPEC(x,\Delta)$, around the infeasible point $x^0$:
\begin{mini!}
  {\substack{d\in\R^{n}}}
  {\transp{\nabla f(x)} d}
  {\label{eq:projlpcc}}
  {}
  \addConstraint{0}{\le c(x) + \transp{\nabla_xc(x)}d}
  \addConstraint{0}{\le \Delta - \infnorm{d}}
  \addConstraint{0}{\le x_1+d_1 \perp x_2+d_2\ge 0.}
\end{mini!}
In order to verify the feasibility of the corresponding branch NLP, we then directly try to solve it using the NLP solver of our choice, in this case {\madnlp}.

\begin{algorithm}
  \caption{An LPEC crossover method for {\MPECs}}
  \label{alg:crossover}
  \begin{algorithmic}[1]
    \Input
    \Desc{$\hat{x}$}{A guess for a feasible point of the \MPEC}.
    \State Compute $\Delta^0 = \alpha_{\Delta_\mathrm{cross}}\max\paren{\infnorm{c(\hat{x})}, \infnorm{\hat{X}_1\hat{x}_2}}$
    \State $k = 0$
    \While{$\Delta^k \le \Delta_{\mathrm{mp}}$}
      \State Solve
      \If{$d^k=\PROJLPEC(\hat{x},\Delta^{k})$ is feasible}
        \State $\hat{\I}_1 = \I_1(\hat{x}+d^k)$
        \State $\hat{\I}_2 = \I_2(\hat{x}+d^k)$
      \Else
        \State $\Delta^{k+1} = \alpha_{\mathrm{cross}}\Delta^k$
      \EndIf
      \State $k = k+1$
    \EndWhile
    \If{$\hat{\I}_1, \hat{\I}_2$ not found}
      \Return \texttt{FAILURE}
    \EndIf
    \State $\gamma^1 = \Delta^k$
    \For{$j=1,\ldots, N_{\mathrm{cb}}$}
      \If{$x^j = \BNLP(\hat{\I}_1,\hat{\I}_2,\hat{x},\gamma^j))$~\eqref{eq:bnlp} is feasible}
        \State $x^* = $ \Cref{alg:mpecopt} applied to $x^j$, $\hat{\I}_1$, and $\hat{\I}_2$.
        \State{\Return{$x^*$}}
      \Else
        \State $\gamma^{j+1} = \alpha_{\gamma} \gamma^j$
      \EndIf
    \EndFor
    \State{\Return \texttt{FAILURE}}
  \end{algorithmic}
\end{algorithm}

\subsection{Solving the {\LPCCs}}
As one of the subproblems that this crossover procedure solves is an {\LPCC}, we now describe the several methods we implement for solving such problems.

The first option implemented in {\solver} is using one of either the relaxation or exact penalty algorithms to solve the {\LPCC}.
These methods provide only an $\epstol$-solution to the {\LPCC} however and as such need their own active-set identification procedure to identify $\I_1$ and $\I_2$.
In our implementation, we apply a ``naive'' active-set identification which compares the magnitudes of $x_1$ and $x_2$ elementwise.

The second method we provide is based on the so called ``Big-M'' reformulation of the LPCC into a  Mixed Integer Linear Program (MILP).
The combinatorial nature of complementarity constraints naturally leads to the possibility of reformulating them using integer variables, with one binary variable per complementarity: $z\in \cbrac{0,1}^{\ncc}$.
This means that global solutions and certificates of infeasibility can be found by specialized solvers for this kind of problem.
Furthermore, when solving $\PROJLPEC(x,\Delta)$, these algorithms can be terminated early, that is without verifying global optimality, once a feasible integer assignment is found.
This is due to \cite[Theorem 4.1, Remark 4.1]{Nurkanovic2025}, which says that only a \textit{feasible} point of the {\LPCC} is needed to identify a feasible branch.
We provide interfaces to use the commercial solver {\gurobi}~\cite{Gurobi} and open source solver {\highs}~\cite{Huangfu2018} which can be used to solve moderately sized instances.

While the {\solver} based algorithms do not come with the global optimization and infeasibility detection guarantees that the MILP solvers have, for large problems, we observe that these algorithms begin to become competitive with even commercial MILP solvers.
For huge scale problems, such as those arising in Security-Constrained Optimal Power Flow (SC-OPF), the {\solver} based {\LPCC} solvers become the only feasible option.
We illustrate such an example below.
For sake of brevity, we omit extensive benchmarks for the crossover implementations, but we observe similar performance as reported in~\cite{Nurkanovic2025}.

\begin{example}[Crossover for huge scale problems]\label{ex:crossover}
  In this example, we demonstrate how using the algorithms implemented in {\solver} to solve {\LPCC} subproblem described above opens the door to calculating extremely accurate solutions to huge scale {\MPCCs}.
  In order to do this, we solve the \texttt{case\_ACTIVSg500} SC-OPF problem with 10 contingencies~\cite{Pacaud2025}.
  This problem has 42,960 variables, 53,733 general nonlinear constraints, and 2,220 complementarities.
  In \Cref{tab:crossover_inaccurate} we demonstrate a log of the iterations needed to cross-over to an accurate solution from an $x^0$ acquired using the relaxation algorithm with $\epstol=10^{-4}$.
  For the projection trust region we take a $\Delta^0 = 10^{-1}$.
  This procedure takes on the order $23.3$ seconds during which the mixed-integer based solvers fail to find an integer feasible solution for the projection {\LPCC}.
  For comparison, in \Cref{tab:crossover_accurate} we demonstrate a log of the iterations needed to cross-over to an accurate solution from an $x^0$ acquired using the relaxation algorithm with $\epstol=10^{-8}$ and a projection trust region of $\Delta^0 = 10^{-4}$.
  This process takes a shorter but analogous time of $17.2$ seconds.
  In both cases we used $\Delta = 10^{-4}$ as the B-stationarity verification initial trust region.
  The relaxation solver achieves the required tolerance in $6.1$ and $13.3$ seconds respectively.

  Both projections successfully find a solution with a zero complementarity residual, that is, they satisfy the upper level complementarities $0 \le x_1 \perp x_2 \ge 0$ exactly.
  We note that as is the case in linear and quadratic programming, the crossover procedure is an expensive one.
  While in some cases attempting the {\LPCC} based projection from a very inaccurate solution can actually be faster than achieving an accuracy of $\epstol = 10^{-8}$ using the regularization based solvers, most of the time this is not the case.
  However, the complementarity accuracy achievable using the crossover method is practically infeasible to achieve using other methods.
    Note that if we try to approximate this complementarity accuracy with a lower $\epstol$ (or even just requiring the complementarity tolerance to be smaller on its own) the relaxation and penalty solvers will often stall or fail to compute steps due to the degeneracy of the problem near the solution.
  In applications where accurate selection of the MPCC branches in large scale problems is critical, we show that using {\solver}, makes the computational cost of finding such solutions feasible in the same order of magnitude as achieving the tolerance of $\epstol=10^{-8}$.
 As such, while we pay a moderate computation price for this accuracy, it still significantly improves accuracy when this is required.
  \begin{table}[t]
    \centering
    \scriptsize
    \begin{tabular}{c|| c c c c c c p{3cm}}
      Iter & \# {\LPCC} & \# BNLP & $|I_{00}|$ & $\Delta$ & $||\textrm{Step}||$ & $||\Delta f||$ & Description\\\hline
      $0$ & $1$ & $1$ & $8$ & $10^{-1}$ & -- & -- & Projection.\\
      $1$ & $1$ & $1$ & $2$ & $10^{-4}$ & $0.54$ & $4.5\times 10^{-6}$ & BNLP step\\
      $2$ & $1$ & $1$ & $2$ & $10^{-4}$ & -- &$-5.0\times 10^{-7}$ & BNLP step rejected\\
      $3$ & $2$ & $0$ & $2$ & $10^{-6}$ & -- & -- & B-stat. verified
    \end{tabular}
    \caption{Crossover iterations for $\epstol = 10^{-4}$ solution in Phase I.}
    \label{tab:crossover_inaccurate}
  \end{table}

    \begin{table}[t]
    \centering
    \scriptsize
    \begin{tabular}{c|| c c c c c c p{3cm}}
      Iter & \# {\LPCC} & \# BNLP & $|I_{00}|$ & $\Delta$ & $||\textrm{Step}||$ & $||\Delta f||$ & Description\\\hline
      $0$ & $1$ & $1$ & $1$ & $10^{-4}$ & -- & -- & Projection.\\
      $1$ & $1$ & $1$ & $1$ & $10^{-4}$ & -- & $-5.6\times 10^{-8}$ & BNLP step rejected\\
      $2$ & $1$ & $1$ & $1$ & $10^{-5}$ & -- &$-1.2\times 10^{-8}$ & BNLP step rejected\\
      $3$ & $1$ & $0$ & $1$ & $10^{-6}$ & -- & -- & B-stat. verified
    \end{tabular}
    \caption{Crossover iterations for $\epstol = 10^{-8}$ solution in Phase I.}
    \label{tab:crossover_accurate}
  \end{table}
\end{example}

\section{Benchmarking}
\label{sec:benchmarks}
In this section we first discuss the software packages and tools used to generate the following numerical results.
We then present several comparative benchmarks using the \MacMPEC~\cite{Leyffer2000} benchmark suite
and three larger benchmarks: (i) security-constrained optimal power flow (SC-OPF) instances~\cite{Pacaud2025},
(ii) a set of problems taken from {\NOSBENCH}~\cite{Nurkanovic2024b},
(iii) and a suite of Quadratic Programs with Complementarity Constraints (QPCCs) arising in nonsmooth model predictive control~\cite{Nurkanovic2026a}.
We use the larger problems in order to better understand the performance of our algorithms, as many of the problems in {\MacMPEC} are small and pose little challenge even to naive solution methods.

In order to aid in reproducibility we collect all the software used to generate the following results in one git repository at: \url{https://github.com/apozharski/ccopt_paper_examples.git}.

\subsection{Problem formulation}
The solvers we have implemented in \solver\ solve \MPCCs\ in the form:
\begin{mini!}[4]
  {\substack{x\in\R^{n}}}
  {f(x)}
  {\protect\label{eq:madmpec_mpec}}
  {}
  \addConstraint{l_g}{\le g(x) \le u_g\protect\label{eq:madmpec_mpec:g}}
  \addConstraint{l_{x_0}}{\le x_0 \le u_{x_0}}
  \addConstraint{l_{x_1}}{\le x_1 \perp x_2 \ge l_{x_2}}
  \addConstraint{u_{x_1}}{\ge x_1}
  \addConstraint{u_{x_2}}{\ge x_2,}
\end{mini!}
with $l_g,u_g\in \R_\infty^{n_g}$, $l_{x_0},u_{x_0}\in \R_\infty^{\nzero}$, $u_{x_1},u_{x_2}\in \R_\infty^{\ncc}$, and $l_{x_1},l_{x_2}\in\R^{\ncc}$.
Finite values for lower bounds are required for the well-posedness of the complementarity constraints, but all other bounds may be set to $\pm\infty$ to represent one sided bounds.
We note that {\solver} does not currently support the solution of mixed-complementarity constraints, which can be written compactly as $a \perp \ell \leq b \leq u$.

\subsection{Solvers}
This section contains a brief description of the solvers we compare against {\solver}.
In cases where options differ for a particular test set, we describe these changes in the corresponding section.

\paragraph{Homotopy.}
A classically successful approach to solving large scale {\MPCCs} has been the direct homotopy method~\cite{Nurkanovic2024b}.
In this method the relaxed NLP:
\begin{mini*}[4]
  {\substack{x\in\R^{n}}}
  {f(x)}
  {}
  {}
  \addConstraint{l_g}{\le g(x) \le u_g\protect\label{eq:madmpec_mpec:g}}
  \addConstraint{l_{x_0}}{\le x_0 \le u_{x_0}}
  \addConstraint{0}{\ge x_1 x_2 - e\tau^k}
  \addConstraint{l_{x_2}}{\le x_1 }
  \addConstraint{u_{x_1}}{\ge x_1}
  \addConstraint{u_{x_2}}{\ge x_2,}
\end{mini*}
is solved to convergence for a given regularization parameter $\tau^k$. Then the regularization
$\tau$ is updated as $\tau^{k+1} = 0.1\tau^k$ and a new NLP is solved. We set $\tau^0 = 1$.
We provide an implementation of this method in the {\solver} package with a limited feature set.
This implementation can use either {\ipopt}~\cite{Waechter2009} or {\madnlp}~\cite{Shin2024} as the underlying NLP solver.
For the algorithmic parameters we take the values as discussed in~\cite{Nurkanovic2024b}, and use the quality function based barrier update rule in both {\ipopt} and {\madnlp}, which have been shown to have good performance in practice.
For all of the benchmarks in this section we use {\ipopt} as the NLP solver and \texttt{HSL MA27}~\cite{Duff1982} as the linear system solver.
The homotopy solver is explicitly referred to as {\ipopt}-Homotopy.

\paragraph{{\lcqpow}~\cite{Hall2024}.}
For the convex QPCC test sets we compare against several versions of the {\lcqpow} algorithm.
Namely we use the two sparsity exploiting versions of the algorithm.
(i) The first of these uses the quadratic programming solver qpOASES~\cite{Ferreau2014}, which in turn relies on the sparse linear solver \texttt{HSL MA57}~\cite{Duff2004}.
(ii) The second uses the sparsity exploiting quadratic programming solver OSQP~\cite{Stellato2020}.
We use the default settings provided in~\cite{Hall2024} for both versions.

\paragraph{{\gurobi}~\cite{Gurobi}.}
In practice, it is popular to solve problems with complementarity constraints using mixed-integer approaches~\cite{Aussel2025,Siddiqui2013}.
When feasible, we evaluate our solver against a MILP solver using the state-of-the-art commercial software package {\gurobi}.
The first approach uses {\gurobi}'s SOS1 constraint functionality.
For this algorithm we use a node limit of 1000, a \texttt{MIPGap} of $10^{-4}$ and a \texttt{FeasibilityTol} of $10^{-6}$.
The second approach, which we use only for the spring mass friction system dataset from {\lcqpow}, uses the mixed-integer reformulation from~\cite[2.2.3]{Hall2024}.
For this case, we use the same options as provided in the LCQPTest repository\footnote{\label{url:lcqptest}\url{https://github.com/apozharski/LCQPTest}}.

\paragraph{Penalty Homotopy.}
For the comparison against previously published results, we also include an {\ipopt} based penalty homotopy algorithm.
The details of this algorithm can be found in~\cite[2.2.3]{Hall2024} and the LCQPTest repository\footnoteref{url:lcqptest}.

\paragraph{{\knitro}~\cite{Byrd2006}.}
Finally, for large scale SC-OPF problems we compare against the commercial NLP solver package {\knitro}.
{\knitro} uses by default the linear solver HSL MA57, and implements the penalty interior-point method described in \cite{Leyffer2006}.

\subsection{Benchmarks}
Now we describe the four benchmarks we use to evaluate the performance of {\solver}.
For each we discuss the results and compare performance to other relevant solvers.

\subsubsection{\MacMPEC\ benchmark}
As our initial verification {\solver} we use the {\MacMPEC} problem set.
The problem set is imported into Julia via the \texttt{AmplNLReader.jl} julia package~\cite{Orban2022}.
While the relative ease of this benchmark has been recently documented~\cite{Nurkanovic2025,Nurkanovic2024b}, it provides a good initial baseline and a proof of concept for any {\MPCC} solver implementation.
We solve 180 problems from the test set, only omitting problems formulated as mixed complementarity problems.

\begin{figure}[t]
  \centering
  \begin{subfigure}{0.49\linewidth}
    \centering
    \includegraphics[width=\textwidth]{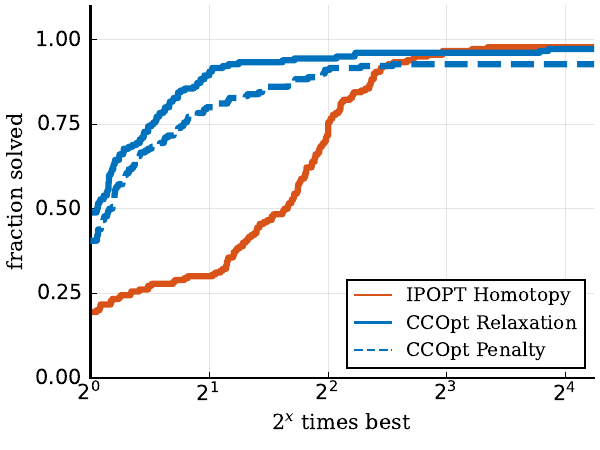}
    \caption{Relative performance on {\MacMPEC} in terms of iterations.}
  \end{subfigure}
  \hfill
  \begin{subfigure}{0.49\linewidth}
    \includegraphics[width=\textwidth]{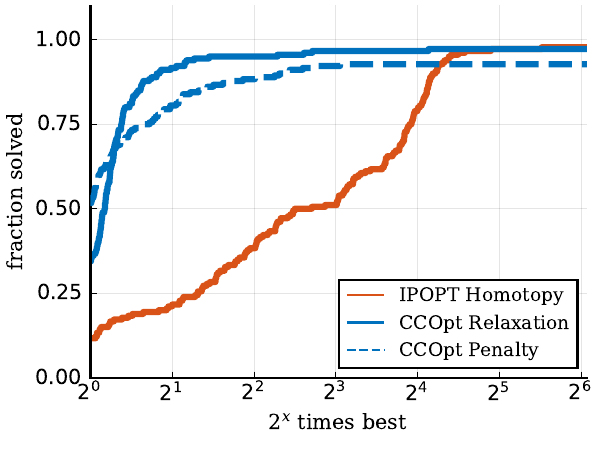}
    \caption{Relative performance on {\MacMPEC} in terms of wall time.}
  \end{subfigure}
  \caption{A comparison of solvers on the {\MacMPEC} benchmark.}
  \label{fig:macmpec}
\end{figure}

The performance plots in \Cref{fig:macmpec} show a wall time comparison as well as an iteration count comparison.
As all the algorithms use an interior-point style method and solve similar linear systems, the iteration count becomes a fair metric to compare solver efficiency.
We note that four of the problems in this test suite are known to be infeasible and are indicated as ``failures'' in our benchmark.These are namely the problems: \texttt{pack-rig2-16}, \texttt{pack-rig2-32}, \texttt{pack-rig2c-16}, and \texttt{pack-rig2c-32}.

When the {\solver}-Relaxation method and the {\ipopt}-Homotopy solvers converge, they find the best known solution in all cases.
{\solver}-Penalty solver occasionally converges to feasible solutions which are marginally worse than the best known solution.
We note the strong performance of both the {\solver}-Relaxation and {\solver}-Penalty algorithms, achieving in many cases a 2 to 20 times speedup in terms of wall time, and a 2 to 8 times reduction in iteration count when compared to {\ipopt}-Homotopy.
We also see a minimal robustness reduction particularly when using {\solver}-Relaxation.
Namely, the \texttt{monteiro} and \texttt{monteiroB} instances suffer from great sensitivity to initialization and fail due to unbounded iterates.
These problems are Bilevel problems and have degenerate (non-complementarity) constraint Jacobians, which lead to unbounded steps in the dual variables without regularization.
We observe similar issues when using {\madnlp} as the NLP solver in the classical homotopy algorithm to solve this problem.
Both problems can, however, be easily solved if a small fixed dual regularization $\delta_c$ is chosen.

It is clear to us from these results that {\solver} exhibits good performance on this test set,
and as such we proceed to test on a suite of large-scale MPCCs coming from relevant applications.

\subsubsection{SC-OPF benchmark}
We solve a set of large-scale SC-OPF problems formulated as {\MPCCs}.
These problems are described in~\cite{Pacaud2025}: here, the SC-OPF computes the optimal dispatch for a set of generators in a transmission network, while maintaining feasibility with relation to a set of contingency scenarios (line tripping).
The problems are implemented using \texttt{ExaModels.jl}~\cite{Shin2024}, for fast function evaluation.
The complementarity constraints model the adjustment of the generators in the post-contingency states, and the PV/PQ switches.
The number of complementarity constraints scale linearly with the number of contingencies $K$, and the larger problems often formulate with more than 10,000 complementarity constraints.
As such they make a good test case for our solver, and help us evaluate the performance on real world large scale {\MPCCs}.

\begin{figure}[t]
  \centering
  \begin{subfigure}{0.49\linewidth}
    \centering
    \includegraphics[width=\textwidth]{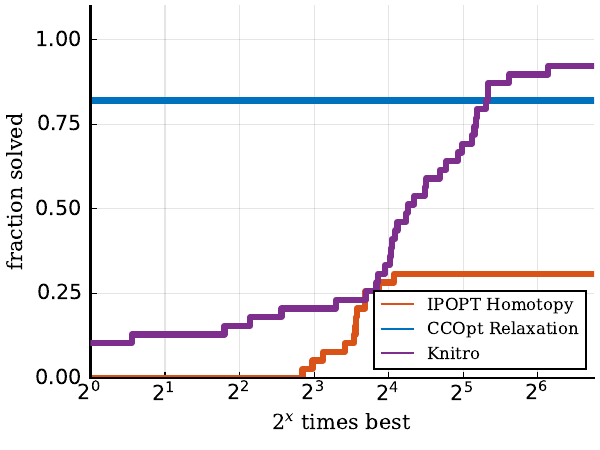}
    \caption{Relative performance on SC-OPF instances in terms of wall-time.}
  \end{subfigure}
  \hfill
  \begin{subfigure}{0.49\linewidth}
    \includegraphics[width=\textwidth]{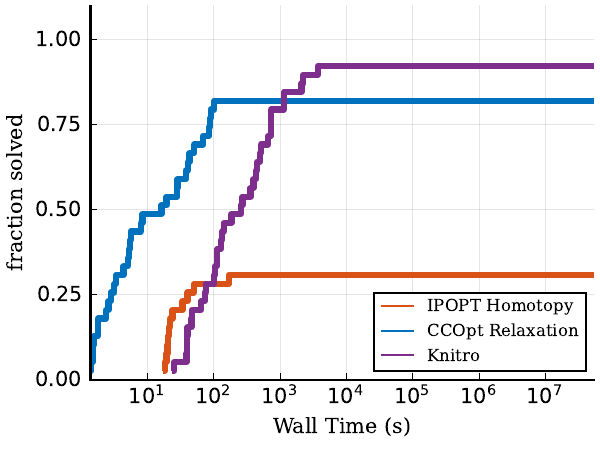}
    \caption{Absolute performance on SC-OPF instances in terms of wall time.}
  \end{subfigure}
  \caption{A comparison of solvers on the SC-OPF instances.}
  \label{fig:scopf}
\end{figure}

We benchmark {\solver} on 40 SC-OPF instances from the GO-competition~\cite{Aravena2023}. We consider three networks with 200, 500 and 2000 buses, and increase the number of contingencies from 10 to 100.
The largest instance has 400k variables and 20k complementarity constraints.
We compare {\solver}-Relaxation with {\ipopt}-Homotopy and with {\knitro}.
All the solvers use the linear system solver \texttt{HSL MA27}~\cite{Duff1982}, and the instances are solved with a tolerance $\epstol = 10^{-8}$.
The two performance profiles are displayed in Figure \ref{fig:scopf}.
We observe that {\ipopt}-Homotopy solves only 12 out of 40 instances, whereas {\solver}-Relaxation method and Knitro solve respectively 32 and 36 instances.
Despite being slightly less robust than Knitro, {\solver}-Relaxation is significantly faster when it converges: the average solution time is 24s for {\solver}-Relaxation, compared to 393s for Knitro.

\subsubsection{\NOSBENCH}
A second application field for {\MPCCs} is the optimal control and simulation of nonsmooth systems~\cite{Lin2022,Nurkanovic2024a}.
The difficulty of such {\MPCCs} when compared to classical benchmarks has been one of the driving forces behind the development of {\solver}.
We use a representative subsample of 90 problems from the {\NOSBENCH}~\cite{Nurkanovic2024b} benchmark suite to evaluate the performance of our solver.
For the exact list of problems we point the reader to the {\NOSBENCH} code repository \url{https://github.com/nosnoc/nosbench}.
As noted in~\cite{Nurkanovic2024b}, the optimal control instances in this test suite are particularly difficult and solving them efficiently remains an open challenge.
In this benchmark all solvers use \texttt{HSL MA27}~\cite{Duff1982} solver for the linear system solves and a tolerance of $\epstol = 10^{-7}$.
We compare against the {\ipopt}-Homotopy algorithm implemented in {\nosnoc} which has been previously shown in~\cite{Nurkanovic2024b} to be the most effective for this suite.
For a more accurate comparison when measuring the performance of the {\nosnoc} solver, we omit the overhead of the Matlab implementation and report only the time spent in the NLP solver.
This is done as the unpredictable nature of the Matlab runtime behavior disproportionally affects {\ipopt}-Homotopy as it returns to the Matlab runtime between each NLP solve.
\begin{figure}[t]
  \centering
  \begin{subfigure}{0.49\linewidth}
    \centering
    \includegraphics[width=\textwidth]{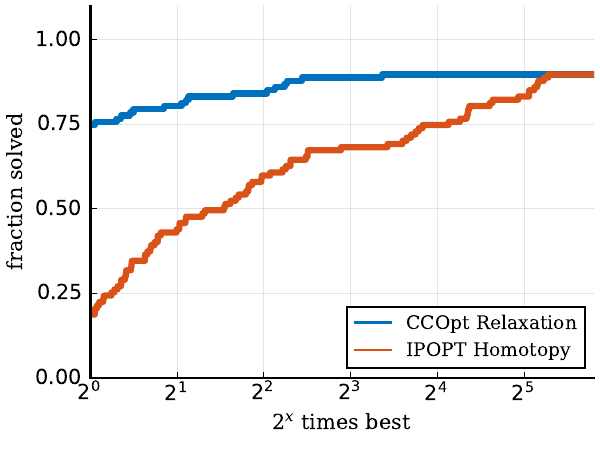}
    \caption{Relative performance on {\NOSBENCH} in terms of wall time.}
  \end{subfigure}
  \hfill
  \begin{subfigure}{0.49\linewidth}
    \includegraphics[width=\textwidth]{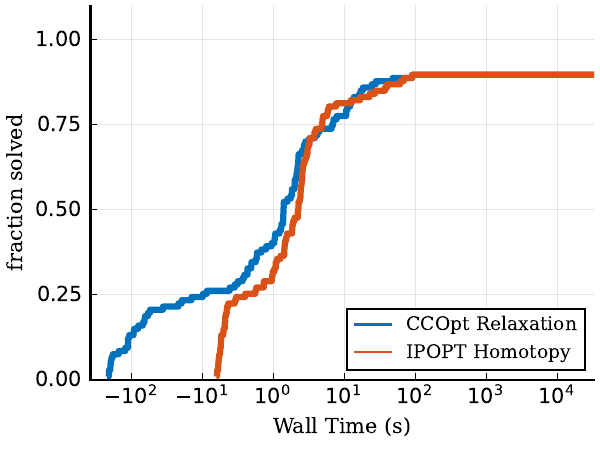}
    \caption{Absolute performance on {\NOSBENCH} in terms of wall time.}
  \end{subfigure}
  \caption{A comparison of solvers on the {\NOSBENCH} benchmark.}
  \label{fig:nosbench}
\end{figure}

In this benchmark we see a major performance uplift over {\ipopt}-Homotopy.
{\solver}-Relaxation obtains a wall time performance increase of 2 to 30 times, and with improved robustness.
We see the largest improvements in performance on simulation problems where locally isolated solutions exist.

\subsubsection{Quadratic programs with complementarity constraints}
It has been recently shown that Quadratic Programs with Complementarity Constraints (QPCCs) seem to be the correct tool for real-time Model Predictive Control (MPC) algorithms~\cite{Nurkanovic2026a}.
For real-time applications, solution performance is critical to allow reduced sampling rates and rapid feedback, which are in turn critical for effective MPC.
As such we evaluate our solver on a test suite of 251 QPCCs which come from the MPC examples in {\nosnoc}~\cite{Nurkanovic2022b}.
We split this problem set into two distinct groups: the first with 124 problems coming from QPCCs whose Hessians have been convexified via the ``mirroring'' technique~\cite{Verschueren2021}, and the second with 127 problems without convexification applied (and thus some may have indefinite Hessians).

We evaluate performance of {\solver} against {\ipopt}-Homotopy in {\nosnoc} as well as the commercial solver {\gurobi} using its SOS1 constraint functionality.
For the problems with convexified Hessians we further compare against the dedicated QPCC solver {\lcqpow}.
We solve these QPCCs to a tolerance of $\epstol = 10^{-6}$ with all solvers.

\begin{figure}[t]
  \centering
  \begin{subfigure}{0.49\linewidth}
    \centering
    \includegraphics[width=\textwidth]{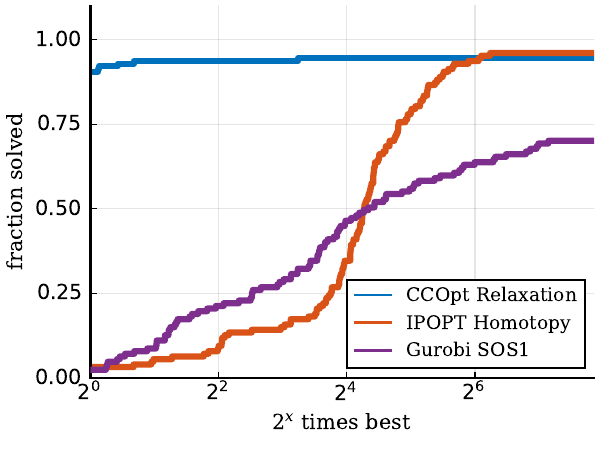}
    \caption{Relative performance on \NOSBENCH\texttt{-QPCC} in terms of wall time.}
  \end{subfigure}
  \hfill
  \begin{subfigure}{0.49\linewidth}
    \includegraphics[width=\textwidth]{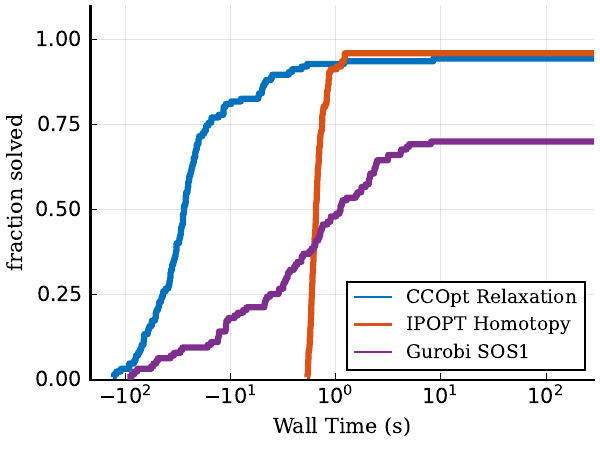}
    \caption{Absolute performance on \NOSBENCH\texttt{-QPCC} in terms of wall time.}
  \end{subfigure}
  \caption{A comparison of solvers on the \NOSBENCH\texttt{-QPCC} benchmark.}
  \label{fig:qpcc}
\end{figure}

\begin{figure}[t]
  \centering
  \begin{subfigure}{0.49\linewidth}
    \centering
    \includegraphics[width=\textwidth]{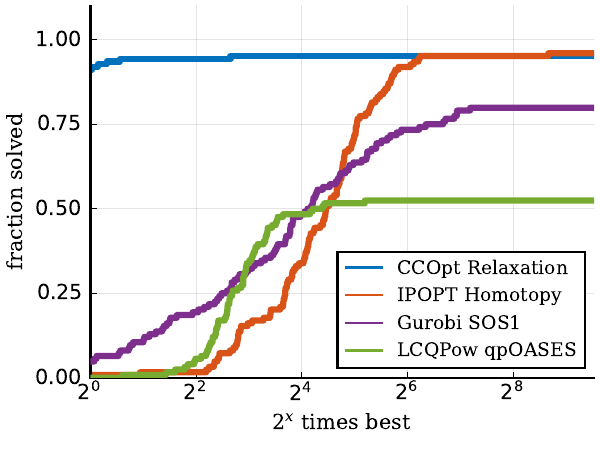}
    \caption{Relative performance on \NOSBENCH\texttt{-QPCC-CONVEX} in terms of wall time.}
  \end{subfigure}
  \hfill
  \begin{subfigure}{0.49\linewidth}
    \includegraphics[width=\textwidth]{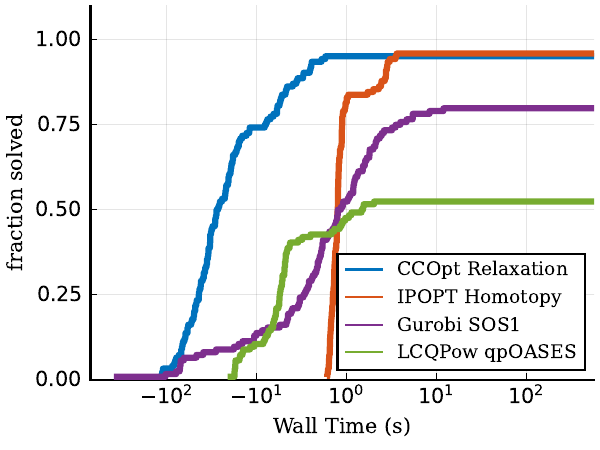}
    \caption{Absolute performance on \NOSBENCH\texttt{-QPCC-CONVEX} in terms of wall time.}
  \end{subfigure}
  \caption{A comparison of solvers on the \NOSBENCH\texttt{-QPCC-CONVEX} benchmark.}
  \label{fig:qpcc_convex}
\end{figure}

The performance profile for the first two benchmarks are pictured in \Cref{fig:qpcc_convex} and \Cref{fig:qpcc} respectively.
In both  we observe that {\solver}-Relaxation nearly fully dominates all other methods with speedups of between 2 and 64 times while being nearly as robust as {\ipopt}-Homotopy.
We emphasize that unlike {\lcqpow}, the only specialization to QPCCs we use in {\solver} is the Q-regularization scheme.
In particular, this means that the Hessian of the Lagrangian and constraint Jacobian are re-evaluated at every new iterate in the algorithm, which for many QPCC solves, constitutes nearly half the run time.
This speedup shows that {\solver}-Relaxation is a significant step towards real time optimal control of nonsmooth systems.

In order to verify the validity of our findings further, we use an existing test set of QPCCs from the moving masses test set in~\cite[Section 5.3]{Hall2024}.
These QPCCs are a formulation of an optimal control problem described by Stewart in~\cite{Stewart1996}.
We evaluate the sparse qpOASES and OSQP versions of {\lcqpow}, {\ipopt}-Homotopy, and {\gurobi}.
For all solvers but {\solver}, we use the same options as those presented in the originating publication.
The implementation of these can be found in the git repository at \url{https://github.com/apozharski/LCQPTest} and in~\cite{Hall2024}.
We observe a one to two order of magnitude improvement over the comparably robust solvers, and a factor of two to ten improvement over the OSQP based version of {\lcqpow}, which is significantly less robust.

\begin{figure}[t]
  \centering
  \begin{subfigure}{0.49\linewidth}
    \centering
    \includegraphics[width=\textwidth]{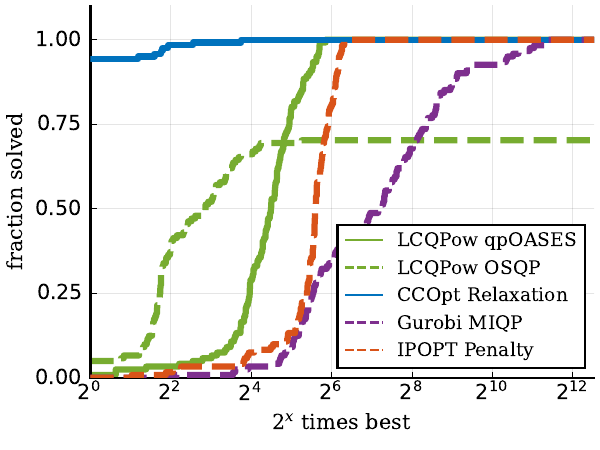}
    \caption{Relative performance on the Moving Masses dataset.}
  \end{subfigure}
  \hfill
  \begin{subfigure}{0.49\linewidth}
    \includegraphics[width=\textwidth]{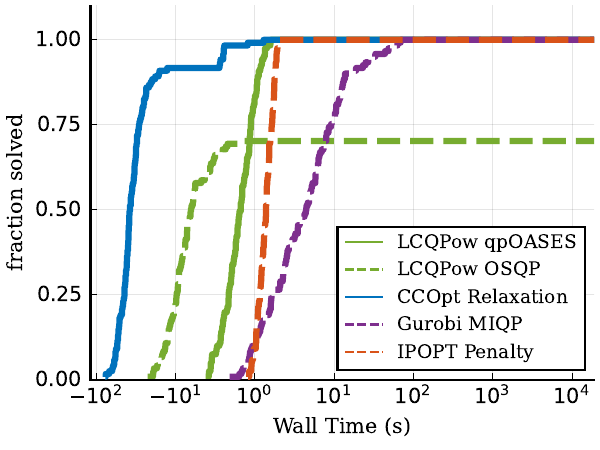}
    \caption{Absolute performance on the Moving Masses dataset.}
  \end{subfigure}
  \caption{A comparison of solvers on the Moving Masses dataset.}
  \label{fig:moving_masses}
\end{figure}


\section{Conclusion and future work}
\label{sec:conclusion}
In this work, we detail the implementation of {\solver}, a high performance solver suite tailored to Mathematical Programs with Complementarity Constraints (\MPCCs).
The suite includes two implementations of interior-point Methods (IPM) for {\MPCCs}, one using a relaxation reformulation, and the other using a penalty based reformulation.
The solvers are built on the IPM infrastructure of \madnlp.
As part of these implementations we describe several algorithmic improvements focusing on several key areas which constitute the major difficulties in solving \MPCCs.
These are namely: ill-conditioning and indefiniteness of the Lagrange Hessian in the KKT matrix, and the selection of relaxation, penalty, and barrier parameters as the algorithms search for a solution.
For the former we propose and implement a Hessian regularization scheme which exploits the known structure of the complementarity constraints to alleviate their indefinite contribution to the Hessian of the Lagrangian.
Additionally, we propose an ``endgame'' algorithm which selectively relaxes the lower bounds of the complementarity variables, and thus avoids having a nearly empty interior for the complementarity constraint reformulation.
This often improves the conditioning of the KKT system.
For the latter, we implement an \MPCC\ tailored optimization-based parameter update scheme which updates simultaneously the barrier and relaxation parameters in the relaxation-based algorithm.
We implement also an adaptive penalty update for the penalty-based algorithm.
Numerically, {\solver}'s performance is often an entire order of magnitude better than existing algorithms, including commercial solvers, on a suite of large scale problems arising from several applications of {\MPCCs} in science and engineering, .
Finally, we implement the recently proposed \texttt{mpecopt} algorithm~\cite{Nurkanovic2025}, allowing for high accuracy solutions and verification of B-stationarity with regards to the found solution.

These algorithms have been extensively benchmarked and good defaults have been selected for all algorithmic parameters, a subset of which can be found in \Cref{sec:solver_options}.
This in concert with the implementation of a compiled shared library interface and an interface to \texttt{CasADi}, allows practitioners to easily use the software using their choice of Julia, Python, Matlab, or C++.
The software package is released under the MIT Open Source license and source code is available at \url{https://github.com/MadNLP/CCOpt.jl}

While extensive benchmarking suggests the implementation of this solver represents a significant step forward in open source solvers for {\MPCCs} there are still several significant avenues for future work.
The first of these is the specialization of this solver and its algorithms to {\LPCCs} and quadratic programs with complementarity constraints.
These subclasses of {\MPCCs} still constitute nonconvex nonlinear programs due to the presence of complementarity constraints, however their structure has exploitable properties.
The second avenue to improve the solver's performance is a tighter integration with the numerical linear algebra routines used to calculate search directions.
One possible improvement is the integration of GPU accelerated linear solvers such as NVIDIA's CUDSS in conjunction with more accurate methods of step calculation such as HyKKT~\cite{Regev2023}.
Further, structure exploiting factorizations similar to classic Ricatti-recursion-based methods can be developed for problems with block-banded KKT systems such as those which appear in the optimal control of nonsmooth systems.

\backmatter

\bmhead{Supplementary information}

Data tables containing benchmark results are provided as supplementary materials.

\bmhead{Acknowledgements}

The authors would like to thank Sungho Shin, Alexis Montoison, and Tim Besard for their prompt code review and merging of work needed for the ahead-of-time compiled library containing {\solver} and {\madnlp}.
Further, thanks is given to Joris Gillis and Joel Andersson for their aid with and review of the {\casadi} interface for {\solver

\section*{Declarations}

\bmhead{Funding}
This research was supported by DFG via projects 504452366 (SPP 2364), 560056112 (robust MPC), 535860958 (ALeSCo) and 525018088 (MAWERO), and by BMWK via 03EN3054B.
\bmhead{Competing interests}
The authors have no relevant financial or non-financial interests to disclose
\bmhead{Code availability}
The code discussed in this paper are available at the links provided above.

\clearpage
\bibliography{syscop}

\appendix
\section{Summary of {\solver} options}\label{sec:solver_options}
In this section we provide a brief reference for the algorithmic options of {\solver}.
A more in depth list can be found in the {\solver} code and documentation: \url{https://github.com/MadNLP/CCOpt.jl}.
\newpage
\begin{landscape}
  \subsection{Relevant \madnlp\ Options}
  \begin{longtable}[c]{l c c c p{5cm}}
    Option Name & Symbol & Default Value & Relevant Section & Description\\\hline
    \texttt{barrier.mu\_init} & $\mu^0$ & $0.1$ & \Cref{sec:ipm} & Initial barrier parameter.\\
    \texttt{tol} & $\epstol$ & $10^{-8}$ & \Cref{sec:madnlpc:term_crit} & Stationarity tolerance.\\
    \texttt{barrier} & -- &\texttt{MonotoneUpdate} & \Cref{sec:madnlpc:homotopy_update} & Type of barrier updated to use.
  \end{longtable}

  \subsection{Relevant \solver-Relaxation Options}
  \begin{longtable}[c]{l c c c p{6cm}}
    \caption{Main relaxation solver options}
    \endfirsthead
    \endhead
    Option Name & Symbol & Default Value & Section & Description\\\hline
    \texttt{q\_regularization} & -- & \texttt{critical\_rho} & \Cref{sec:ipm} & Q regularization scheme to use.\\
    \texttt{critical\_rho\_factor} & $\alpha_B$ & $0.9999$ & \Cref{sec:madnlpc:hessian_regularization} & Factor of critical multiplier used.\\
    \texttt{min\_eig\_value} & $\lambda_{\mathrm{min}}$ & $10^{-8}$ & \Cref{sec:madnlpc:hessian_regularization} & Minimum eigenvalue of complentarity contribution in Lagrangian Hessian.\\
    \texttt{relaxation\_update} & -- & \texttt{RolloffRelaxationUpdate} & \Cref{sec:madnlpc:homotopy_update} & Update rule to use for $\tau$.\\
    \texttt{endgame\_strategy} & -- & \texttt{RelaxLBEndgameStrategy} & \Cref{sec:madnlpc:endgame} & Endgame algorithm employed.\\
    \texttt{endgame\_threshold} & $\epsilon_{\mathrm{endgame}}$ & $10^{-6}$ & \Cref{sec:madnlpc:endgame} & KKT error threshold at which endgame algorithm is triggered.\\
    \texttt{center\_complementarities} & -- & \texttt{true} & \Cref{sec:initialization} & Whether the complementarity variables should be centered.\\
    \texttt{centering\_factor} & $k_{\mathrm{cen}}$ & $0.5$ & \Cref{sec:initialization} & How far along the $x_1 = x_2$ line to place the initial centered complementarities.\\
    \texttt{mu\_thresh} & -- & $5\times 10^{-6}$ & -- & The barrier value below which the kkt regularization described in \cite{Raghunathan2005} is enabled.\\
  \end{longtable}

  \begin{longtable}[c]{l c c c p{6cm}}
    \caption{\texttt{ProportionalRelaxationUpdate} options}
    \endfirsthead
    \endhead
    Option Name & Symbol & Default Value & Section & Description\\\hline
    \texttt{sigma\_mu\_ratio} & $\alpha_\tau$ & $1.0$ & \Cref{sec:madnlpc:homotopy_update} & Proportional factor for $\tau$.\\
    \texttt{sigma\_mu\_exp} & $\beta_\tau$ & $1.0$ & \Cref{sec:madnlpc:homotopy_update} & Exponential factor for $\tau$.\\
    \texttt{sigma\_min} & -- & $10^{-8}$ & \Cref{sec:madnlpc:homotopy_update} & Minimum value for $\tau$.
  \end{longtable}

  \begin{longtable}[c]{l c c c p{6cm}}
    \caption{\texttt{RolloffRelaxationUpdate} options}
    \endfirsthead
    \endhead
    Option Name & Symbol & Default Value & Section & Description\\\hline
    \texttt{rolloff\_slope} & $a$ & $2.0$ & \Cref{sec:madnlpc:homotopy_update} & Defines slope of $\tau(\mu)$ as $\mu\rightarrow 0$\\
    \texttt{rolloff\_point} & $b$ & $10^{-6}$ & \Cref{sec:madnlpc:homotopy_update} & Defines (partially) at what $\mu$ we begin to reduce $\tau$.\\
    \texttt{rolloff\_max} & $c$ & $1.0$ & \Cref{sec:madnlpc:homotopy_update} & Maximum value for $\tau$.\\
    \texttt{sigma\_min} & -- & $10^{-8}$ & \Cref{sec:madnlpc:homotopy_update} & Minimum value for $\tau$.
  \end{longtable}

  \begin{longtable}[c]{l c c c p{6cm}}
    \caption{\texttt{LOQORelaxationUpdate} options}
    \endfirsthead
    \endhead
    Option Name & Symbol & Default Value & Section & Description\\\hline
    \texttt{gamma} & $\gamma_{\tau}$ & $2.0$ & \Cref{sec:madnlpc:homotopy_update} & Factor for LOQO rule.\\
    \texttt{r} & -- & $10^{-8}$ & \Cref{sec:madnlpc:homotopy_update} & LOQO rule step size.\\
    \texttt{sigma\_min} & -- & $10^{-8}$ & \Cref{sec:madnlpc:homotopy_update} & Minimum value for $\tau$.
  \end{longtable}
  \newpage
  \begin{longtable}[c]{l c c c p{6cm}}
    \caption{\texttt{RelaxLBEndgameStrategy} options}
    \endfirsthead
    \endhead
    Option Name & Symbol & Default Value & Section & Description\\\hline
    \texttt{delta\_max} & $\delta_{\mathrm{max}}$ & $10^{-4}$ & \Cref{sec:madnlpc:endgame} & Maximum allowed lower bound relaxation.\\
    \texttt{tau} & $\xi$ & $10^{-8}$ & \Cref{sec:madnlpc:endgame} & Exponent applied to $\infty$ norm of residual to identify whether engame should be triggered.
  \end{longtable}

  \subsection{Relevant \solver-Penalty Options}
  \begin{longtable}[c]{l c c c p{6cm}}
    \caption{Main penalty solver options}
    \endfirsthead
    \endhead
    Option Name & Symbol & Default Value & Section & Description\\\hline
    \texttt{rho\_0} & $\rho_0$ & $1.0$ & \Cref{sec:penalty} & Initial penalty.\\
    \texttt{rho\_max} & -- & $10^{10}$ & \Cref{sec:penalty} & Maximum penalty.\\
    \texttt{rho\_growth\_rate} & $\alpha_\rho$ & $10$ & \Cref{sec:ell1:parameter_update} & Increase factor for penalty.\\
    \texttt{gamma} & $\gamma$ & $0.4$ & \Cref{sec:ell1:parameter_update} & Exponent in expression used for triggering dynamic penalty update.\\
    \texttt{q\_regularization} & -- & \texttt{critical\_rho} & \Cref{sec:ipm} & Q regularization scheme to use.\\
    \texttt{critical\_rho\_factor} & $\alpha_B$ & $0.99$ & \Cref{sec:madnlpc:hessian_regularization} & Factor of critical penalty to use.\\
    \texttt{min\_eig\_value} & $\lambda_{\mathrm{min}}$ & $10^{-8}$ & \Cref{sec:madnlpc:hessian_regularization} & Minimum eigenvalue of complentarity contribution in Lagrangian Hessian.\\
    \texttt{comp\_history\_length} & $h$ & $10$ & \Cref{sec:ell1:parameter_update} & Length of history kept for complementarity values.\\
    \texttt{eta\_dynamic\_update} & $\eta_{\mathrm{pen}}$ & $0.99$ & \Cref{sec:ell1:parameter_update} & Sufficient decrease parameter for dynamic penalty updates.
  \end{longtable}
\end{landscape}

\end{document}